\documentclass{amsart}
\usepackage{amssymb}
\usepackage{amsfonts}
\usepackage{amsbsy}

\usepackage{latexsym}
\usepackage[mathscr]{eucal}
\usepackage{verbatim}
\usepackage{color}
\theoremstyle{plain}
\newtheorem{theorem}{Theorem}[section]

\newtheorem{corollary}[theorem]{Corollary}
\newtheorem{lemma}[theorem]{Lemma}
\newtheorem{proposition}[theorem]{Proposition}
\newtheorem{remark}[theorem]{Remark}
\newtheorem{definition}[theorem]{Definition}
\newtheorem{example}[theorem]{Example}

\newcommand{\A}{\mathcal{A}}
\newcommand{\B}{\mathcal{B}}

\newcommand{\R}{\mathbb{R}}

\newcommand{\Z}{\mathcal{Z}}
\newcommand{\K}{\mathcal{K}}

\newcommand{\eps}{\epsilon}

\DeclareMathOperator{\TA} {\mathcal{T}(\A)}
\DeclareMathOperator{\Aff} {Aff(\TA)}
\DeclareMathOperator{\Affk} {Aff(K)}
\DeclareMathOperator{\LAff} {LAff(\TA)}
\DeclareMathOperator{\LAffk} {LAff(K)}
\DeclareMathOperator{\LAffs} {LAff_\sigma(\TA)}
\DeclareMathOperator{\LAffsk} {LAff_\sigma(K)}
\DeclareMathOperator{\LAffb}{LAff_b(\TA)}
\DeclareMathOperator{\LAffq} {LAff(\mathcal{QT}(\A))}

\DeclareMathOperator{\Ext} {\partial_{e}(\TA)}
\DeclareMathOperator{\Extk} {\partial_{e}(K)}

\DeclareMathOperator{\Her} {Her}
\DeclareMathOperator{\tr} {Tr}

\DeclareMathOperator{\St} {\A\otimes \K}
\DeclareMathOperator{\M}{\mathcal M(\mathcal A \otimes \mathcal K)}
\newcommand{\Ma}{\mathcal M(\mathcal A)}
\newcommand{\Mb}{\mathcal M(\mathcal B)}

\DeclareMathOperator{\Mul}{\mathcal M}
\DeclareMathOperator{\Ped}{Ped(\A)}
\DeclareMathOperator{\co}{co}
\newcommand{\Imin}{I_{\mathrm{min}}}
\newcommand{\Icon}{I_{\mathrm{cont}}}
\newcommand{\Ifin}{I_{\mathrm{fin}}}
\newcommand{\Ib}{I_{\mathrm{b}}}
\newcommand{\+}{\overset{\cdot}{+}}
\newcommand{\Sc}{\mathscr S}

\def\sideremark#1{\ifvmode\leavevmode\fi\vadjust{\vbox to0pt{\vss
\hbox to 0pt{\hskip\hsize\hskip1em
\vbox{\hsize2cm\tiny\raggedright\pretolerance10000
\noindent#1\hfill}\hss}\vbox to8pt{\vfil}\vss}}}

\newcommand{\be}{\begin{equation}\label}
\newcommand{\ee}{\end{equation}}
\newcommand{\bq}{\begin{equation*}}
\newcommand{\eq}{\end{equation*}}
\newcommand{\ba}{\begin{align*}}
\newcommand{\ea}{\end{align*}}
\newcommand{\bp}{\begin{proof}}
\newcommand{\ep}{\end{proof}}
\newcommand{\bL}{\begin{lemma}\label}
\newcommand{\eL}{\end{lemma}}
\newcommand{\bP}{\begin{proposition}\label}
\newcommand{\eP}{\end{proposition}}
\newcommand{\bC}{\begin{corollary}\label}
\newcommand{\eC}{\end{corollary}}
\newcommand{\bT}{\begin{theorem}\label}
\newcommand{\eT}{\end{theorem}}
\newcommand{\bR}{\begin{remark}\label}
\newcommand{\eR}{\end{remark}}
\newcommand{\bD}{\begin{definition}\label}
\newcommand{\eD}{\end{definition}}
\newcommand{\bE}{\begin{example}\label}
\newcommand{\eE}{\end{example}}
\numberwithin{equation}{section}

\thanks{This work was partially supported by the Simons Foundation (grant No 245660 to Victor Kaftal and grant No  281966 to Shuang Zhang).}

\author{Victor Kaftal}

\address{Department of Mathematics\\
University of Cincinnati\\
P. O. Box 210025\\
Cincinnati, OH\\
45221-0025\\
USA}

\email{kaftalv@ucmail.uc.edu}

\author{P.  W. Ng}

\address{Department of Mathematics\\
University of Louisiana\\
217 Maxim D. Doucet Hall\\
P. O. Box 3568\\
Lafayette, Louisiana\\
70504-3568\\
USA}

\email{png@louisiana.edu}

\author{Shuang Zhang}

\address{Department of Mathematics\\
University of Cincinnati\\
P.O. Box 210025\\
Cincinnati, OH\\
45221-0025\\
USA}

\email{zhangs@ucmail.uc.edu}

\date{}


\begin{document}
\title[ ]{Purely infinite corona algebras}

\begin{abstract}

Let $\A$ be a simple, $\sigma-$unital, non-unital  C*-algebra, with metrizable tracial simplex $\TA$, projection-surjectivity and injectivity, and strict comparison of positive elements by traces. Then the following are equivalent: 
\begin{enumerate}
\item [(i)] $\A$ has quasicontinuous scale;
\item [(ii)] $\Ma$ has strict comparison of positive elements by traces;
\item [(iii)] $\Ma/\A$ is purely infinite;
\item [(iii$'$)] $\Ma/\Imin$ is purely infinite;
\item [(iv)] $\Ma$ has finitely many ideals;
\item [(v)] $\Imin=\Ifin$.\\
If furthermore $M_n(\A)$ has projection-surjectivity and injectivity for every $n$, then the above conditions are equivalent to:
\item [(vi)] $V(\Ma)$ has finitely many order ideals.
\end{enumerate}
Quasicontinuity of the scale is a notion introduced by Kucerovsky and Perera that extends both the property of having finitely many extremal traces and of having continuous scale.  Projection-surjectivity and injectivity permit to identify projections in $\Ma\setminus \A$ with lower semicontinuous affine functions on $\TA$. $\Imin$ is the smallest ideal of $\Ma$ properly containing $\A$,  and   $\Ifin$ is the ideal of of $\Ma$ generated by the positive elements with evaluation functions finite over the extremal boundary of $\TA$
\end{abstract}

\maketitle

\section{Introduction}
In the study of multiplier algebras of C*-algebras, an important role is played by the associated corona algebras.  In the case of the algebra $\K$ of compact operators on a separable Hilbert space $H$, $\Mul(\K)= B(H)$ and the corona algebra $\Mul(\K)/ \K$ is the Calkin algebra which is well known to be both simple and purely infinite. 

Perhaps one reason for the success of the Brown--Douglas--Fillmore theory 
(\cite{BDF1}, \cite{BDF2}) is that in their context, the multiplier algebra $\Mul(\K)$ 
and the corona algebra $\Mul(\K)/ \K$ have particularly nice structure. For example, the BDF--Voiculescu result which, roughly says
that every essential extension is absorbing, would not be true if the corona
algebra $\Mul(\K)/ \K$ were not simple (\cite{ArvesonDuke}, \cite{BDF1},
\cite{Voiculescu}).  Thus it is natural that structure and comparison theory
for multiplier algebras and corona algebras are indispensable
(though often implicitly present) for operator theory and extension theory in this general  context. For example, it is by now clear that, in the classical theory of absorbing extensions, ``nice" extension theory corresponds to ``nice" corona algebra structure (e.g., \cite{ElliottKucerovsky}, \cite{Kasparov1}, \cite{Kasparov2},
\cite{Lin1},  \cite{LinStableUniqueness}, \cite{Paschke}, \cite{Zhang1}).

If $\A$ is simple, $\sigma$-unital but not unital, and non-elementary C*-algebra, Lin showed in \cite {LinContScale} and \cite {LinSimple} that $\Ma/\A$ is simple if and only if $\Ma/\A$ is simple and purely infinite, if and only if $\A$ has continuous scale. Simple continuous scale algebras are one of the most interestingl classes in nonstable generalizations of BDF Theory.

However, the continuity of the scale  is not necessary for $\Ma/\A$ to be purely infinite. Indeed, Kucerovsky and Perera \cite  {KucPer} identified a weaker condition, which they called quasicontinuity of the scale,  which is necessary and sufficient for  $\Ma/\A$ to be purely infinite in the case that is $\A$ simple, separable, non-unital, with real rank zero, stable rank one, strict comparison of positive elements by quasitraces, and finitely many  infinite extremal quasitraces.  
Their result was then extended by Perera,  Ng, and Kucerovsky \cite [Theorem A] {KucNgPer} to simple, separable C*-algebras that are the stabilization of a unital C*-algebra and are either exact and $\Z$-stable or are AH-algebras with slow dimension growth.  Furthermore, they showed that quasicontinuity of the scale  is sufficient for the same class of C*-algebras in order for $\Ma$ to have finitely many ideals  and it is also necessary in case $\A$ is exact.

The notion of quasicontinuity of the scale (see Definition \ref {D:quasicontinuous scale}) can be extended to any C*-algebra with non-empty tracial simplex. As the terminology suggests, for  C*-algebras with non-empty tracial simplex, continuous scales are quasicontinuous. It is also easy to see that if a simple C*-algebra is the stabilization of a unital algebra (equivalently, is stable and contains a non-zero projection), then the scale is quasicontinuous if and only if the algebra has only finitely many extremal traces. 

The main goal of this paper is to present the connection between quasicontinuity of the scale of $\A$, pure infiniteness of the corona algebra $\Ma/\A$, and other properties of $\Ma$ that are essentially connected with the first two.

We start by showing in Proposition \ref {P:strict comp implies prop inf} that  if $\Ma$ has strict comparison of positive elements by traces (see Definition \ref {D:str comp}),  then $\Ma/\A$ is purely infinite.

In \cite [Theorem 6.6 ] {KNZCompPos} we have proven that if $\A$ is simple, $\sigma$-unital, has strict comparison of positive elements by traces, and has quasicontinuous scale, then $\Ma$ has strict comparison of positive elements by traces. As a consequence, quasicontinuity of the scale and strict comparison of positive elements by traces for $M_3(\A)$ imply that $\Ma/\A$ is purely infinite, (implications (i) $\Rightarrow$ (ii) $\Rightarrow$ (iii) of Theorem \ref{T:main}).

Also, using the  technique developed  in \cite [Theorem 6.6 ] {KNZCompPos}, we prove in Corollary \ref{C: structure ideals} that if $\A$ is a simple, $\sigma$-unital, non-unital    C*-algebra, with quasicontinuous scale and strict comparison of positive elements by traces, then $\Ma$ has finitely many ideals, (implication (i) $\Rightarrow$ (iv) of Theorem \ref{T:main}).

To prove that quasicontinuity of the scale is also necessary for these properties  (items (ii), (iii), and (v)  of Theorem \ref{T:main}), we need to assume additional regularity conditions for $\A$. Notice that the real rank zero algebras considered in \cite {KucPer}, have plenty of projections, while the Jiang-Su algebra $\Z$ considered in \cite {KucNgPer} does not contain non-trivial projections.. The common thread is that  both classes of algebras are projection-surjective and injective  in the sense of Definitions \ref  {D:surject and inject}, that is, the evaluation map $\hat P: \TA\ni \tau \mapsto \tau(P)$  establishes a 1:1 and onto correspondence between the equivalence classes of projections in $\Ma\setminus \A$ and strictly positive lower semicontinuous  affine functions on the tracial simplex $\TA$ that are pointwise limits of an increasing sequence of continuous affine functions and that are complemented under the scale  $ \Sc:=\widehat{1_{\Ma}}$.  A  large class of simple C*-algebras are projection-surjective and projection-injective (see Section \ref {S: proj inj surj}). For simple, separable, stable algebras, projection-surjectivity and injectivity implies strict comparison of positive elements of the algebra (Theorem \ref {T: proj surj inj implies strict comp}); and  both are equivalent if the algebra has also stable rank one and contains a nonzero projection (Corollary \ref {C:surjectivity}). 
A key part in this study is also played by three distinguished ideals of $\Ma$, $\Imin$, $\Icon$, and $\Ifin$. $\Imin$ (see (\ref {e:def Imin})) is the smallest ideal of $\Ma$ that properly contains $\A$, and was studied in \cite {LinIdealsAF}, \cite{LinZhang}, \cite{LinContScale}, \cite {LinSimple}, \cite {PereraIdeals}, and more recently in \cite {KNZMin}. If $\A$ is separable or if $\A$ has strict comparison, then $\A\ne \Imin$ and $\Imin/\A$ is purely infinite and simple (\cite [Corollary 3.15, Theorem 4.8] {KNZMin}). $\Icon$ is the ideal generated by the elements with continuous evaluation functions and it coincides with $\Imin$ if $\A$ has strict comparison. Without strict comparison we have an example where the two are different (\cite[Theorem 5.6, Theorem 7.8] {KNZMin}).  Another important ideal is $\Ifin$ (see Definition \ref {D:ideals}) which was called the {finite ideal} in \cite {PereraIdeals} and is the ideal generated by the elements with evaluation function finite on the extremal boundary $\Ext$ of the tracial simplex $\TA$.

For C*-algebras with strict comparison and projection-surjectivity and injectivity, we show that strict comparison on $\Ma$, pure infiniteness of $\Ma/\A$, and the finiteness of the ideal lattice of $\Ma$ are equivalent to quasicontinuity of the scale, by proving that they are equivalent to $\Imin=\Ifin$.

The paper is organized as follows: In Section 2, we present some background  material on the tracial simplex and on ideals of $\Ma$. In Section 3, we present some technical lemmas on lower semicontinuous affine functions.  In Section 4, we introduce quasicontinuity of the scale and its relation to the finiteness of the ideal lattice of $\Ma$. In Section 5, we introduce and discuss the notion of projection-surjectivity and projection-injectivity and how they are used to obtain results on ideals in multiplier algerbas. In Section 6, we prove our main result (Theorem \ref {T:main})  linking all the various equivalent conditions. Then, under additional assumptions on the tracial simplex, we derive further  properties of the ideal lattice of $\Ma$; in particular, the existence of infinite chain of principal ideals.

\section{Notations and preliminaries}
\subsection{The tracial simplex and strict comparison}
Given a simple, $\sigma$-unital, (possibly unital)  C*-algebra $\A$   and   a nonzero positive element $e$ in the Pedersen ideal $\Ped$ of $\A$, denote by $\TA$  the collection of the (norm) lower semicontinuous  densely defined tracial weights $\tau$ on $\A_+$, that are normalized on $e$. Explicitly, a trace $\tau\in \TA$ is an additive and homogeneous map from $\A_+$   into $[0, \infty]$ (a weight); satisfies the trace condition $\tau(xx^*)=\tau(x^*x)$ for all $x\in \A$; is densely defined (also called densely finite, or semifinite), i.e., the positive cone $\{x\in \A_+\mid  \tau(x)< \infty\}$  is dense in $\A_+$; satisfies the lower semicontinuity condition $\tau(x)\le \varliminf \tau(x_n)$ for $x, x_n\in \A_+$ with $\|x_n-x\|\to 0$, or equivalently,  $\tau(x)= \lim_n \tau(x_n)$ for $0\le  x_n\uparrow x$ in norm; and is normalized on $e$, i.e., $\tau(e)=1$. 
We will mostly assume that $\TA\ne \emptyset$ and hence that $\A$ is stably finite. 

When equipped with the topology of pointwise convergence on $\Ped$,  $\TA$ is a Choquet simplex (e.g., see \cite [Proposition 3.4]{TikuisisToms} and \cite {ElliottRobertSantiago}). In particular, $\TA$ is a compact convex  subset of a locally convex linear topological Hausdorff space, {\it compact convex  space} for short. 
 The collection of the extreme points of $\TA$ is denoted by $\Ext$ and is called the {\it extremal boundary} of $\TA$. For simplicity's sake we call the elements of $\TA $ (resp., $\Ext$) traces (resp., extremal traces.) Tracial simplexes $\TA$ arising from different nonzero positive elements in $\Ped$ are homeomorphic; so we will not reference explicitly which element $e$ is used for the normalization. A trace $\tau$ on $\A$ extends naturally  to a trace on $\St$ (explicitly to the trace $\tau\otimes\tr $), and so we can identify  $\mathcal{T}(\St)$ with $\TA$. 
By the work of F. Combes   \cite [Proposition 4.1, Proposition 4.4] {Combes} and Ortega, R{\o}rdam, and Thiel \cite[Proposition 5.2]{OrtegaRordamThiel},
every trace $\tau \in \TA$ has a unique extension to a lower semicontinuous (i.e., normal) tracial weight (trace for short)  on the enveloping von Neumann algebra $\A^{**}$, and hence to a trace on the multiplier algebra $\Ma$. We will still denote that extension by $\tau$. For more details, see  \cite {TikuisisToms}, \cite {ElliottRobertSantiago} and also  \cite {KNZCompPos} and  \cite {KNZComm}.

Although we will use the following notions mainly for the case when $K$ is a Choquet simplex, it is customary (and more convenient) to formulate them for  compact  convex spaces.
\bD{D:Aff}
Given a compact convex  space $K$,
\item [(i)]$\Affk$ denotes the Banach space of the continuous real-valued affine functions on $K$ with the uniform norm;
\item [(ii)]$\LAffk$ denotes the collection of the lower semicontinuous affine functions on $K$ with values in $\mathbb R\cup\{+\infty\}$;
\item [(iii)] $\Affk_{+}$ (resp., $\LAffk_{+}$) denotes the cone of the positive functions (i.e., $f(x)\ge 0$ for all $x \in K$) in $\Affk$ (resp., in $\LAffk$).
\item [(iv)] $\Affk_{++}$ (resp., $\LAffk_{++}$) denotes the cone of the strictly positive functions (i.e., $f(x)>0$ for all $x \in K$) in $\Affk$ (resp., in $\LAffk$).
\item [(v)] $\LAffsk$, (resp., $\LAffsk_{+}$, $\LAffsk_{++}$) denotes the collection of functions in $\LAffk$ (resp., $\LAffk_{+}$, $\LAffk_{++}$) that are the increasing pointwise limit of a \emph{sequence} of functions in $\Affk$.
\item [(vi)] Given $\Sc\in \LAffsk_{++}$, an element $f\in \LAffsk_{++}$ is said to be complemented under $\Sc$ if there is a $g\in \LAffsk_{++}\sqcup \{0\}$ such that $f+g =\Sc$.
 \eD
For every $ A\in \Ma_+$, $\hat A$ denotes the evaluation map and $\widehat{[A]}$ the dimension function of $A$: 
\begin{align}\label{e: eval map} \TA&\ni \tau\mapsto \hat A(\tau):=\tau(A) \in [0, \infty];\\
\label{e: dim map} \TA&\ni \tau\mapsto \widehat{[A]}(\tau):=d_\tau(A)= \lim_n \tau(A^{1/n}) \in [0, \infty].\end{align}
As shown  in \cite[Remark 5.3]{OrtegaRordamThiel}, \be{e:range proj} d_\tau(A)=\tau(R_A)\quad\text{where }R_A\in \A^{**}\text{ is the range projection of }A.\ee 
It is well known that   $\hat A\in \LAffs_{++}$ and $ \widehat{[A]} \in \LAffs_{++}$ for every $A\ne 0$.
 In particular the  {\it scale} $ \Sc$ of $\A$ is defined as $ \Sc := \widehat{1_{\Ma}}$  and  thus
\be{e:Sc} \Sc \in \LAffs_{++}.\ee

We will also use frequently the following well known facts. If $A, B \in \Ma_+$, and $\tau\in \TA$ then
\begin{align}\label {e: ineq->ineq}  A\leq  B~ &\Rightarrow ~ \hat A(\tau) \le \hat B(\tau)\\
\label{e: suub->ineq}
A\preceq B~ &\Rightarrow ~ d_\tau(A)\le d_\tau(B) \\
  \intertext{where  $``\preceq"$ denotes Cuntz subequivalence, that is the existence of a sequence $X_n\in \Ma$ such that $\|X_nBX_n^*- A\|\to 0.$}
AB=0~ &\Rightarrow ~ d_\tau(A+B)=d_\tau(A)+d_\tau(A)\label{e: dim sum}\\
\tau(A)&\le \|A\|d_\tau (A) \label{e:ineq 1}\\
 d_\tau((A-\delta)_+)&< \frac{1}{\delta}\hat A(\tau)~\,~\forall  ~\delta >0 \label{e:ineq 2}
  \intertext{ and by \cite[Lemma 2.4 (iii)]{KNZCompPos}, }
 d_\tau((A+B-\delta_1-\delta_2)_+)&\le d_\tau((A-\delta_1)_+)+d_\tau((B-\delta_2)_+)\quad \forall \delta_1\ge0,\,\delta_2\ge0.  \label{e:ineq 3}
\end{align}
 By the definition of the topology on $\TA$, if $a\in \Ped$, then $\hat a\in \Aff$. Notice that $\widehat{[a]}$ is not necessarily continuous. We will use the fact that  $\widehat{[a]}$ is bounded:
\bL{L:bounded dim} Let $\A$ be a C*-algebra with non empty tracial simplex $\TA$ and let $a\in \Ped_+$. Then $\underset{\tau\in \TA}{\sup} d_{\tau}(a) < \infty.$ 
\eL
\bp
Since $a\le \sum_{j=1}^na_j$ for some $n\in \mathbb N$ elements $a_j\in \A_+$, with {\it local unit}, i.e.,  such that $b_ja_j=a_j$ for some $b_j\in \A_+$, and since $ d_\tau (a)\le \sum_{j=1}^n d_\tau(a_j)$,  it is enough to verify the claim for an $a\in \A_+$ that has a local unit  $b\in \A_+$ (i.e.,  $ba=a$).
Assume without loss of generality that $\|a\|=1$. Since $a$ and $b$ commute, we can identify them as continuous functions on a compact space $X$, i.e., $a=a(x)$ and $b=b(x)$. Then for all $x\in X$ such that $a(x)\ne 0$ we have $b(x)=1$ and hence $f_{\frac{1}{2}}(b(x))=1$, where for every $\eps>0$, the function $f_\epsilon(t)$ is definite as follows:
\be {e:feps}f_\epsilon(t):=\begin{cases}0&\text{for }t\in [0, \epsilon]\\
\frac{t-\eps}{\eps}&\text{for }t\in (\eps, 2\eps]\\
1&\text{for }t\in (2\epsilon, \infty).
\end{cases} \ee

 But then $f_{\frac{1}{2}}(b)a=a$ and $R_a\le f_{\frac{1}{2}}(b)$. Thus  $f_{\frac{1}{2}}(b)$ is also a local unit for $a$ and since itself belongs to $\Ped_+$ as $f_{\frac{1}{4}}(b)f_{\frac{1}{2}}(b)=f_{\frac{1}{2}}(b)$,  its evaluation function $\widehat{f_{\frac{1}{2}}(b)}$ is continuous on $\TA$. Thus
$d_\tau(a)= \tau(R_a)\le \tau\big(  f_{\frac{1}{2}}(b)\big)\quad \forall \tau\in \TA$  and hence $\underset{\tau\in \TA}{\sup} d_{\tau}(a) < \infty.$ 

\ep 
The same result was obtained in \cite[Lemma 1.6] {KucNgPer}  under the additional conditions that $\A$ is the stabilization of a unital simple exact algebra with strict comparison.

The notions of strict comparison has played an important role in the theory of C*-algebras especially after \cite {BlackadarComp}. 
\bD{D:str comp}
Let $\A$ be a simple C*-algebra with $\TA\ne \emptyset$. Then we say that 
\item [(i)] $\A$ has strict comparison of positive elements by traces if  
$a\preceq b$ whenever $a, b\in \A_+$ and $d_\tau(a)< d_\tau(b)$ for all those $\tau\in \TA$ for which $d_\tau(b)< \infty$.
\item [(ii)] $\Ma$  has strict comparison of positive elements by traces if $A\preceq B$ whenever
$A, B\in \Ma_+$,  $A$ belongs to the ideal $I(B)$ generated by $B$, and  $d_\tau(A)< d_\tau(B)$ for all those $\tau\in \TA$ for which $d_\tau(B)< \infty$. 
\eD
Strict comparison is often defined in terms of 2-quasitraces. In \cite [Theorem 2.9] {KNZComm} we proved  that if a unital simple C*-algebra of real rank zero and stable rank one has strict comparison of positive elements by traces (equivalently, of projections, due to real rank zero)  then all 2-quasitraces are traces.  Recently it was shown  the same conclusion holds without the real rank zero and stable rank one hypotheses (\cite[Theorem 3.6] {NgRobert}).

Notice that in (ii), the condition that $A\in I(B)$ (which is obviously necessary for $A\preceq B$) does not follow in general from the condition that $d_\tau(A)< d_\tau(B)$ for all those $\tau\in \TA$ for which $d_\tau(B)< \infty$. Indeed if there is an element $B\in \A_+$ with  $d_\tau(B)= \infty$ for all $\tau\in \TA$ (and this  is certainly the case when  $\A$ is stable) then  the condition $d_\tau(A)< d_\tau(B)$ for all those $\tau\in \TA$ for which $d_\tau(B)< \infty$ is trivially satisfied for every $A\in \Ma_+\setminus \A$  and yet $A\not \preceq B$.

\subsection{Ideals and traces} \label{subS:id}

We first recall the following well-known result.

\bL{L:ideals}
Let $\B$ be a C*-algebra and let $A, T\in \B_+$. Then $A\in I(T)$ (the principal ideal generated by $T$) if and only if for every $\epsilon >0$ there is some  $m\in \mathbb N$ such that $(A-\eps)_+\preceq \bigoplus_{k=1}^m T$ in $M_m(\B)$. In particular, if $P$ is a projection, then  $P\in I(T)$  if and only if there is an $m\in \mathbb N$ such that  $P \preceq \bigoplus_{k=1}^m T$ in $M_m(\B)$.
\eL

We will focus on the ideals of the multiplier algebra $\Ma$ of a simple, non-unital  C*-algebra $\A$. 
The ideal \be {e:def Imin}\Imin:=\bigcap\{ \mathcal J\lhd \Ma\mid \A\subsetneq \mathcal J\}\ee 
is called the minimal ideal of $\Ma$ and $A\subset \Imin$. We do not know in general whether $\A\ne \Imin$ although by \cite[Corollary 3.15, Proposition 5.4, Theorem 5.6]{KNZMin}  this conclusion holds when $\A$ is simple, $\sigma$-unital, non-unital, non-elementary C*-algebra and with any of the following properties:
\begin{enumerate} 
\item $\A$ is separable;
\item the Cuntz semigroup of $\A$ is order separable;
\item $\A$ has the (SP) property and its dimension semigroup $D(\A)$ of Murray-von Neumann equivalence classes of projections is order separable;
\item $\A$ has strict comparison of positive elements by traces.
\end{enumerate}
The conclusion $\A\ne \Imin$ holds also if $\A$ has continuous scale (in particular, if $\A$ is purely infinite), because then $\Ma/\A$ is simple  (\cite [Theorem  2.8] {LinSimple}) and hence $\Imin=\Ma$. 

Following Lin's approach, but not using his notations  (\cite {LinIdeals}), one can characterize $\Imin$ in terms of approximate identities of $\A$. 
Given any approximate identity $\{e_n\}$, which henceforth we will always assume to satisfy the condition $e_{n+1}e_n=e_n$,  the  ideal $\Imin$ is shown (see \cite {LinIdeals} and \cite {KNZMin}) to coincide with the  norm closure of the linear span of  \begin{multline}\notag K_o(\{e_n\}):=\{X\in \Ma_+\mid \forall~0\ne a\in \A_+ ~\exists~N\in \mathbb N \\ \ni ~m>n \ge N\Rightarrow (e_m-e_n)X(e_m-e_n)\preceq a\}.\end{multline}  

When $\A$  is a $\sigma$-unital, non-unital, non-elementary, C*-algebra with non-empty tracial simplex, then another natural ideal is $\Icon$.
\bD{D:Icon}\cite [Definition 5.1, Proposition 5.2]{KNZMin} $\Icon$ is 
the norm closure of the linear span of
$ K_c:= \{X\in \Ma_+\mid \hat X\in \Aff\}.$
\eD
An immediate consequence of the definition (see \cite [Proposition 5.2]{KNZMin}) is that if $0\ne X\in \Ma_+$ and $0\ne P\in \Ma$ is a projection, then
\begin{align}\label{e:Icon}X\in (\Icon)_+&\quad \text{if and only if}\quad \widehat {(X-\delta)_+}\in \Aff ~\forall ~\delta>0;\\
\label{e:cont proj} P\in (\Icon)_+&\quad \text{if and only if}\quad \hat P\in \Aff_{++}. 
\end{align}
There are simple, separable, non-unital  C*-algebras where $\Imin\ne \Icon$ (\cite [Theorem 7.8] {KNZMin}), however $\Imin= \Icon$ when $\A$ has strict comparison of positive elements \cite [Theorem 5.6]{KNZMin}. 

It is well known that every trace $\tau$ gives rise to a (not necessarily proper) ideal 
$I_\tau$ which is the norm closure of the linear span of the hereditary cone
$$\{X\in \Ma_+\mid \tau(X)< \infty\}.$$
As a consequence, if $0\ne X\in \Ma_+$ and $0\ne P\in \Ma$ is a projection, then
\begin{align}\label{e:Itau}X\in (I_\tau)_+&\quad \text{if and only if}\quad  \tau((X-\delta)_+)< \infty ~\forall ~\delta>0;\\
\label{e:proj in Itau} P\in (I_\tau)_+&\quad \text{if and only if}\quad \tau(P)< \infty.
\end{align}
 
In this paper the following ideals will play an important role.
{\bD{D:ideals} Let $\A$  be a $\sigma$-unital, non-unital, non-elementary, C*-algebra with non-empty tracial simplex. 
\item[(i)]$\Ifin:= \cap _{\tau\in \partial_{e}(\TA)}I_\tau$;
\item[(ii)] $\Ib:= \cap _{\tau\in \TA}I_\tau$.
\eD
Perera introduced in a different way the ideal $\Ifin$, which he called called the {\it finite ideal}, for $\sigma$-unital, non-unital, non-elementary C*-algebras with real rank zero, stable rank 1, and weakly unperforated $K_0$ group (\cite [Proposition 6.1]{PereraIdeals}).
The following inclusion is obvious.
\be{e:inclusion}\Icon\subset \Ib\subset \Ifin.\ee
Also an immediate consequence of the definition and of (\ref {e:Itau}), (\ref{e:proj in Itau}) is: 
\bL{L: char Ifin, Ib} Let $\A$  be a $\sigma$-unital, non-unital, non-elementary, C*-algebra with non-empty tracial simplex. 

\item [(i)] $(\Ifin)_+= \{X\in \Ma_+\mid\tau((X-\delta)_+)< \infty  ~\forall \delta >0, \tau\in \Ext \}$;
\item [(ii)] $(\Ib)_+= \{X\in \Ma_+\mid\tau((X-\delta)_+)< \infty  ~\forall \delta >0, \tau\in \TA \}$.\\
In particular, if $P$ is a projection, then
\item [(iii)] $P\in \Ifin~\Leftrightarrow~\hat P(\tau)< \infty $ for all $\tau \in \Ext$.
\eL
To  explain the notation of $\Ib$, we need to make the following elementary observation. 

\bL{L: bdd LAff} Let $K$ be a compact convex space and let $f\in \LAffk_{++}$. Then
\item [(i)] $\underset{x\in \Extk}{\sup }f(x) = \underset{x\in K}{\sup }f(x)$
\item [(ii)] $\underset{x\in K}{\sup }f(x)< \infty$ if and only if $f(x)< \infty$ for all $x\in K$. 
\eL
\bp \item [(i)] It is obvious that $\underset{x\in \Extk}{\sup }f(x)= \underset{x\in co\Extk}{\sup }f(x)$. Then the conclusion follows from the density of $co\Extk$ in $K$ (the Krein-Millman theorem) and the lower semicontinuity of $f$.
\item [(ii)] The necessity being trivial, assume that $\underset{x\in K}{\sup }f(x)= \infty$ and choose a sequence $x_k\in \Extk$ for which $f(x _k)\ge 2^k$ for all $k$. If  $f(x _k)=\infty$ for some $k$, then we are done, thus assume that  $f(x _k)< \infty$ for all $k$. Let $\mu_k$ be the Dirac measure on $x_k$ and $\mu:= \sum_{k=1}^\infty \frac{\mu_k}{2^k}$. Then $\mu$ is probability measure on $\Extk$. Let $x\in K$ be the corresponding element, i.e., $$g(x)= \int_{\Extk}g(y) d\mu(y) =  \sum_{k=1}^\infty \frac{g(x_k)}{2^k}$$ for all $g\in \Affk$ and hence also  $f(x)= \sum_{k=1}^\infty \frac{f(x_k)}{2^k}=\infty.$
\ep
The argument in (ii) is similar to the one in \cite [Lemma 4.4] {PereraIdeals}.

\bC{C: within Ib}  Let $\A$  be a $\sigma$-unital, non-unital, non-elementary, C*-algebra with non-empty tracial simplex, and 
let $0\ne X\in \Ma_+$ and $0\ne P\in \Ma$ be a projection. Then
\item [(i)] $X\in \Ib$ if and only if $\underset{\tau\in \Ext}{sup}\tau ((X-\delta)_+)< \infty ~\forall \delta >0$;
\item [(ii)] $P\in \Ib$ if and only if $\underset{\tau\in \Ext}{sup}\tau(P)< \infty .$
\eC

\bC{C: fin bdry} 
Let $\A$ be a simple, $\sigma$-unital, non-unital  C*-algebra, with nonempty $\TA$, and with $|\Ext|< \infty$. Then
$\Icon=\Ib= \Ifin$.
\eC

Finally we list  our  notations for order ideals.
If $\B$ is a C*algebra, denote by $V(\B)$ the semigroup of Murray von-Neumann equivalence classes of projections in $M_\infty(\B)$, where $[P]+[Q]:=[P\oplus Q]$ for $P,\,Q\in M_\infty(\B)$. 
An order ideal $H$ of $V(\B)$ is a hereditary sub-semigroup of  $V(\B)$. When $[P]\in V(\B)$, denote the principal order ideal generated by $[P]$ by 
\be{e:principal} I([P]):=\{[R]\in  V(\B)\mid [R]\le n[P]\quad\text{for some }n\in \mathbb N\}.\ee
The connection between principal ideals and principal order ideals generated by a projection of $\B$ is an immediate consequence of Lemma \ref {L:ideals} 
\bL{L: id vs oid} If $P, Q$ are projections in $\B$, then the following are  equivalent: 
\item [(i)] $I(P)\subsetneq I(Q)$;
\item  [(ii)] $\exists\, n$ such that $[P]\le n[Q]$, $\not \exists \,n$ such that $[Q]\le n[P]$;
\item  [(iii)] $I([P])\subsetneq I([Q])$.
\eL

\section{Preliminaries on lower semicontinuous affine functions}
Our paper makes use of some technical results on lower semicontinuous affine functions on Choquet simplexes. We collect them in this section. We start by listing for convenience of reference some results that will be used  throughout the paper.

\bP{P: sup aff} Let $K$ be a compact convex  metrizable space.
\item [(i)]  $\LAffk_{++}=\LAffsk_{++}$ ( \cite [Lemma 4.2] {TikuisisToms}, see also comments before Proposition 4.10 in \cite {PereraIdeals}). In particular, for every $f\in \LAffs_{++}$ there is a decomposition $f=\sum_{j=1}^\infty f_j$ (pointwise convergence) with $f_j\in \Aff_{++}$.
\item [(ii)]  \cite [Choquet Theorem, pg 14] {Phelps} For every $x\in K$ there exists a probability measure $\mu$ on $\Extk$ such that $f(x)= \int_{\Extk}f(t)d\mu(t)$ for all $f\in  \LAffk_{+}$. 
\item [(iii)]
If $f, \,g\in \LAffk_{+}$ and $f(x) \ge g(x)$ (resp.,  $f(x) > g(x)$), (resp., $f(x) = g(x)$)  for all $x\in \Extk$ then  $f(x) \ge g(x)$, (resp., $f(x) > g(x)$), (resp., $f(x) = g(x)$),  for all $x\in K$. \item If furthermore $K$ is a Choquet simplex, then the measure in (ii) is unique  \cite [Choquet Theorem, pg 60] {Phelps}.
\eP
\bT{T:GoodearlHB}\cite[Theorem 11.14, Corollary 11.15]{Goodearl}
Let $K$ be a Choquet simplex, $X \subseteq \partial_e K$ a compact subset of the extremal boundary of $K$,
$f : K \rightarrow \{ -\infty \} \cup \R$ an upper semicontinuous
convex function,
$h : K \rightarrow \R \cup \{ \infty \}$ a lower semicontinuous 
concave function, and $g_0: X\rightarrow \mathbb R$ and continuous function, such that  $f \leq h$ and
$f|_X \leq g_0 \leq h|_X.$
Then there exists a function $g \in \Affk$ such that 
\item [(i)] $f \leq g \leq h$, and
\item [(ii)] $g|_X = g_0$.\\
In particular, every function $g_0 \in C(X, \R)$ has an extension $g \in \Affk$  such that $\|g\|=\|g_o\|$.
 \eT

The following is an elementary observation which we will use in a number of occasions.

\bL{L:facts} Assume that $g= G+F$ where $G$ and $F$ are finite and lower semicontinuous functions on a compact set $K$ and that there is a sequence $K\ni x_n\to x\in K$ such that $g(x_n)\to g(x)$. Then $G(x_n)\to G(x)$.
In particular, if $G,F\in \LAffk_{++}$ and $g:=G+F\in \Affk$, then $G,F\in \Affk_{++}$.
\eL

Now recall that  if $K$ is a simplex, then the complementary face $F'$ of a face $F$  is the union of all the faces of $K$ that are disjoint from $F$. 
A face $F$ of $K$ is said to be {\it split} if $K=F\+F'$ where $\+$ denotes the direct convex sum.  By  \cite [Theorem 11.28]{Goodearl},
\be{e: closed face  is split}
\text{if $K$ is a Choquet simplex then every closed face is split.} 
\ee

It is elementary and most likely well known that if $K$ is a Choquet simplex and $F$ is a split face, then every pair of affine nonnegative functions $f$ on $F$ and $g$ on $F'$  has unique extension to an affine function on $K$.
 We will need to use this fact and some refinements of it, collected in the following lemma. 
\bL{L:direct sum}
 Let $K$ be a Choquet simplex and $F$ a split face.   Assume that $f$ and $g$ are affine nonnegative extended real valued function on $F$ and $F'$ respectively and let $f\overset{\cdot}{+}g$ be the function defined on $K$ as follows: if $k=tx+(1-t)y$ for some  $x\in F$, $y\in F'$, and $0\le t\le 1$ then  
 $$(f\overset{\cdot}{+}g)(k):= \begin{cases}tf(x)+(1-t)g(y)&0< t<1\\
g(y)&t=0\\ f(x)&t=1.\end{cases}$$
 Then 
 \item [(i)] $f\overset{\cdot}{+}g$ is the unique affine function that agrees with $f$ on $F$ and with $g$ on $F'$.\\
Assume  henceforth that
\item [(ii)] Assume that  $F$ is closed,  $f\in \mathrm{LAff}(F)_+$, $ g\in \LAffk_{+}$, and $f(x)\le g(x)$ for all $x\in F$. Then $f\overset{\cdot}{+}g\mid_{F'}\in \LAffk_{+}$ and $ f\overset{\cdot}{+}g\mid_{F'}\le  g$.
\item [(iii)] Assume that  $F$ and $F'$  are closed,  $f\in \mathrm{LAff}(F)_+$, and  $g\in\mathrm{LAff}(F')_+$. Then $ f\overset{\cdot}{+}g\in \LAffk_+$.
 \eL
 \bp
\item [(i)]  Recall that a decomposition $k=tx+(1-t)y$ is unique but for the case when $k\in F$ and then $t=1$, $x=k$, and $y$ is arbitrary  or $k\in F'$ and then $t=0$, $y=k$, and $x$ is arbitrary. Therefore the  function $f\overset{\cdot}{+}g$ is well defined.  Also, the definition given can be simplified by the convention that $0\cdot \infty =0$.  A lengthy straightforward computation 
 shows that $f\overset{\cdot}{+}g$ is indeed affine and that it is the unique affine function that agrees with $f$ on $F$ and with $g$ on $F'$. 
 \item [(ii)]  
 Assume that $k_\lambda \to k$ for some net $k_\lambda\in K$. Let $k_\lambda= t_\lambda x_\lambda+(1-t_\lambda)y_\lambda$ for some $x_\lambda\in F$, $y_\lambda\in F'$ and $t_\lambda\in [0,1]$, and let $k= \alpha x_o+ (1-\alpha)y_o$ for some $x_o\in F$, $y_o\in F'$ and $\alpha\in [0,1]$. By passing if necessary to a subnet, we can assume that
 $(f\overset{\cdot}{+}g)( k_\lambda)$ converges to $\varliminf_\lambda (f\overset{\cdot}{+}g)( k_\lambda)$. 
 Then \be{e:(f+g)(klambda)}
 (f\overset{\cdot}{+}g)(k_\lambda)= t_\lambda f(x_\lambda)+ (1-t_\lambda) g(y_\lambda).
\ee
 By the compactness of $[0,1]$, $F$, and $K$, and by passing if necessary to subnets of  subnets,  which will not affect neither $\varliminf_\lambda (f\overset{\cdot}{+}g)( k_\lambda)$ nor $(f\overset{\cdot}{+}g)( k)$, we can assume that $t_\lambda\to t$, $x_\lambda\to x$, and $y_\lambda\to \beta x'+ (1-\beta)y$ for some $x,x'\in F$, $y\in F'$ and $0\le \beta\le 1$. (Notice that if $F'$ is closed, then $x'=0$). Then 
 \be{e:k} k= tx+ \beta(1-t) x'+ (1-\beta)(1-t)y.\ee We leave to the reader the simpler case when $t=\beta=0$, i.e., $k=y$, and thus assume that $t$ and $\beta$ don't both vanish. Set  
 \be{e: defx''} x'':= \frac{t}{t+ \beta(1-t)}x+ \frac{\beta(1-t)}{t+ \beta(1-t)}x'.\ee Then $x''\in F$ and
 $$k= (t+ \beta(1-t))x''  + (1-\beta)(1-t)y$$
 is the decomposition of $k$ in $F\overset{\cdot}{+}F'$.
 Then  
 \begin{alignat*}{2} (f\overset{\cdot}{+}g)(k)&= (t+ \beta(1-t))f(x'') + (1-\beta)(1-t)g(y)&\text{(by definition of $f\overset{\cdot}{+}g$)}\\
 &=tf(x) + \beta(1-t)f(x')+ (1-\beta)(1-t)g(y)\quad &\text{(by (\ref {e: defx''}) as $f$ is affine)}\\
 &\le tf(x) + \beta(1-t)g(x')+ (1-\beta)(1-t)g(y)&\text{(as $f\le  g$ on $F$)}\\
 &= tf(x)+ (1-t)\big ( g(\beta x'+ (1-\beta)y)\big)&\text{($ g$ is affine)}\\
 &\le \varliminf_\lambda t_\lambda f(x_\lambda) + \varliminf_\lambda(1-t_\lambda) g(y_\lambda)&\text{($f$ and $ g$ are lsc)}\\
 &\le \varliminf_\lambda\big(t_\lambda f(x_\lambda) + (1-t_\lambda)g(y_\lambda)\big)\\
&= \varliminf_\lambda (f\overset{\cdot}{+}g)(k_\lambda)&\text{(by definition of $f\overset{\cdot}{+}g$).}
 \end{alignat*}
 \item [(iii)] Follows from the same proof as in (ii).
 \ep
Notice that if both $F$ and $F'$ are closed and $f$ and $g$ are continuous, then the same computation shows that $f\overset{\cdot}{+}g$ is continuous, which of course is well known (e.g., see \cite[Corollary 11.23]{Goodearl}).

Lemma \ref {L:direct sum} provides a generalization of \cite [Proposition 4.10, Corollaries 4.11-13]{PereraIdeals}  to the case when $F$ is closed but not necessarily finite dimensional, and without requiring the metrizability of $K$.

\bC{C:complementation}
Let $K$ be a metrizable Choquet simplex,  $h\in \LAffk_{++}$, $F\subset K$ a closed face such that $h\mid_ F= \infty$. Then
\item [(i)]  $f\overset{\cdot}{+}\frac{h}{2}\in \mathrm{LAff}(K)_{++}$ is complemented under $h$ for every $f\in \mathrm{LAff}(F)_{++}$.
\item [(ii)] For every $f\in \mathrm{LAff}(F)_{++}$  and  $g\in \mathrm{Aff}(K)_{++}$ such that $f\le g\mid_ F$ and $g(x) < h(x)$ for all $x$, then $f\overset{\cdot}{+}g$ is complemented under $h$. In particular, if $\sup f<  \min h$, then for every $\sup f< \gamma< \min h$,  $f\overset{\cdot}{+}\gamma$ is complemented under $h$, where $\gamma$ denotes the constant function $g(x)=\gamma$.
\item [(iii)] If also $F'$ is closed, $f\in\LAffk_{++}$, and $f$ is continuous on  $F'$, then $f$ is complemented under $nh$ for some $n\in\mathbb N$.
\eC
\bp
\item [(i)] $f\overset{\cdot}{+}\frac{h}{2}\in \mathrm{LAff}(K)_{++}$ by Lemma \ref {L:direct sum} as $f\le h\mid_F$. Moreover 
$$f\overset{\cdot}{+}\frac{h}{2}+\frac{h}{2}=\begin{cases} f(x)+\infty= h(x)&x\in F\\
\frac{h(x)}{2}+\frac{h(x)}{2}=h(x)&x\in F'.\end{cases}$$ and hence $f\overset{\cdot}{+}\frac{h}{2}+\frac{h}{2}=h$. Obviously, $\frac{h}{2}\in \mathrm{LAff}(K)_{++}$.
\item [(ii)]
$f\overset{\cdot}{+}g\in \mathrm{LAff}(K)_{++}$ by Lemma \ref {L:direct sum} and $h-g\in  \mathrm{LAff}(K)_{++}$ 
$$\big(f\overset{\cdot}{+}g+ h-g\big)(x)= \begin{cases} f(x)+ h(x)-g(x)=\infty = h(x)& x\in F\\
g(x) + h(x)-g(x) = h(x) & x\in F'.\end{cases}$$
Therefore $f\overset{\cdot}{+}g+ h-g = h$  and $h-g\in \mathrm{LAff}(K)_{++}$ which concludes the proof.
\item [(iii)] $h$ has a strictly positive minumum on the compact set $F'$, hence we can find $n\in \mathbb N$ such that $ f(x)< n h(x)$ for all $x\in F'$ and  and let $g:= nh\mid_{F'}- f\mid_{F'}$. Then $h\mid _F\+ g\in \LAffk_{++}$ by Lemma \ref {L:direct sum} and reasoning as above, $f+h\mid _F\+ g= nh$.
\ep

\bC{C:increasing sequences} Let $K$ be a Choquet simplex and $\{x_j\}_1^\infty \subset \Extk$ be a sequence with distinct terms. Then for every nondecreasing  sequence  of scalars $0< \alpha_j<\infty$ there is a function $f\in\LAffsk_{++}$ such that 
$f(x_j)=\alpha_j$ for all $j$. Moreover $\alpha_1\le f\le \sup_j\alpha_j.$
\eC
    \bp
Starting with $g_o= \alpha _1>0$ we construct an increasing sequence of functions $g_k\in \Affk_{++}$ such that $g_k(x_j)= \alpha_j$ for all $1\le j\le k$ and $g_k(x)\le \alpha_k$ for all $x\in K$.

Assuming the construction up to $k-1$ for some $k\ge 1$, set $X_k:=\{x_1, x_2, \cdots, x_k\}$ and 
$g_{k,o}(x_j):=\begin{cases} g_{k-1}(x_j)&1\le j<k\\
\alpha_k& j=k
\end{cases}$.
Then $g_{k-1}\le \alpha_{k-1}\le \alpha_k$ and hence $g_{k-1}\mid_{X_k}\le g_{k,o}\le \alpha_k$. Let $g_k\in \Affk$ be the extension of $g_{k,o}$ to $K$ for which
$g_{k-1}\le g_{k}\le \alpha_k$ that is given by Theorem \ref {T:GoodearlHB}. Then  $f:= \lim_k g_k$ satisfies the desired properties.
\ep
  
Next we present two technical constructions of lower semicontinuous functions that will be needed in the study of principal ideals in $\Ma$. 
 \bL{L: under discont}
Let $K$ be a compact metrizable space  and $g$ be a non-negative, finite,  lower semicontinuous function on $K$  that is not continuous at some point $x_o\in K$. Then there is a decomposition $g= G+F$ into the sum of lower-semicontinuous non-negative functions  $G$ and $F$ which are both discontinuous at $x_o$ but for which there is a sequence $y_k\to x_o$ such that $G(y_k)\to G(x_o)$ but $g(y_k)\not \to g(x_o)$.

If furthermore $K$ is a compact  convex metrizable space and $g\in \LAffk_+$ (resp., $g\in \LAffk_{++}$),  then we can choose $G, F$ to be in $\LAffk_{+}$ (resp., in $\LAffk_{++}$).
\eL
  
\bp
Since $g$ is lower semicontinuous and $K$ is metrizable, by Proposition \ref {P: sup aff} we can decompose it into a sum $g=\sum_{k=1}^\infty g_k$ of functions $g_i\in C(K, \mathbb R)_{+}$ (resp., $g_i\in C(K, \mathbb R)_{++}$ if $g$ is strictly positive.) 
Since $g$ is  not continuous at $x_o$, there is a sequence $x_j\to x_o$ and a number $\beta$ such that  $g(x_j)>\beta> g(x_o)$ for all $j$. Let $\delta:= \frac{\beta - g(x_o)}{3}$. 
We construct inductively  two sequences of positive integers $m_j\le n_j< m_{j+1}$ starting with $m_1=1$ and two strictly increasing sequences of integers $s_k$ and $t_k$ such that if we set $G_k:= \sum_{j=1}^k\sum_{i=m_j}^{n_j}g_i$,  we have for all integers $k\ge 1$
\begin{enumerate}
\item [(i)] $G_k(x_{s_k})> G_k(x_o) + \delta$;
\item [(ii)] $|G_k(x_j)- G_k(x_o)| <\frac{\delta}{k}$ for all $j\ge t_k;$
\item [(iii)] $\sum_{i=m_{k+1}}^\infty g_i(x_{t_k})< \frac{\delta}{k}$ and $\sum_{i=m_{k+1}}^\infty g_i(x_o)< \frac{\delta}{k}.$
\end{enumerate}
We start the induction by setting $m_1=1$ and $s_1=1$. Since $g(x_{s_1})>\beta$, choose an integer $n_1\ge 1$  such that $\sum_{i=1}^{n_1}g_i(x_{s_1})> \beta$. Thus  
$$G_1(x_{s_1}) > \beta =g(x_o)+ 3\delta \ge  G_1(x_o)+ 3\delta>G_1(x_o)+ \delta ,$$
thus satisfying condition (i). By the continuity of $G_1$ we can find an index $t_1$ for which (ii) is satisfied. By the convergence of the series $\sum_1^\infty g_i(x)$ for every $x$, choose $m_2> n_1$ so to satisfy (iii). Thus conditions (i)-(iii) are satisfied for $k=1$.

Next assume the construction up to some integer $k$ and notice that this includes the existence of $m_{k+1}>n_k$ that satisfies (iii).  By the continuity of $\sum_{i=1}^{m_{k+1}-1}g_i$ and of $G_k(x)$, choose an integer $s_{k+1}>s_k$ such that for all $j\ge s_{k+1}$ we have
\be{e:<delta} \Big |\sum_{i=1}^{m_{k+1}-1}g_i(x_j)-\sum_{i=1}^{m_{k+1}-1}g_i(x_o)\Big| < \delta~\text{and} ~ |G_k(x_j)-G_k(x_o)|< \delta.\ee
Since $g(x_{s_{k+1}}) >\beta$, choose an integer $n_{k+1}\ge m_{k+1}$ so that
\be{e:> beta}\sum_{i=1}^{n_{k+1}}g_i(x_{s_{k+1}})>\beta.\ee

Then
\begin{align*}
G_{k+1}&(x_{s_{k+1}})- G_{k+1}(x_o)\\
&=G_{k}(x_{s_{k+1}})- G_{k}(x_o) + \sum_{i=m_{k+1}}^{n_{k+1}}g_i(x_{s_{k+1}})-  \sum_{i=m_{k+1}}^{n_{k+1}}g_i(x_o)\\
&= G_{k}(x_{s_{k+1}})- G_{k}(x_o)+ \sum_{i=1}^{n_{k+1}}g_i(x_{s_{k+1}})- \sum_{i=1}^{n_{k+1}}g_i(x_o)\\
& \phantom{abcdef}- \sum_{i=1}^{m_{k+1}-1}g_i(x_{s_{k+1}})+ \sum_{i=1}^{m_{k+1}-1}g_i(x_o)   \\
&>-\delta + \beta - g(x_o) -\delta= \delta
\end{align*}
where \begin{alignat*}{2}
&|G_{k}(x_{s_{k+1}})- G_{k}(x_o)|<\delta &\text{(by (\ref{e:<delta}))}\\
& \sum_{i=1}^{n_{k+1}}g_i(x_{s_{k+1}})>\beta &\text{(by (\ref {e:> beta}))}\\
 &\sum_{i=1}^{n_{k+1}}g_i(x_o) \le g(x_o) &\text{(by the definition of $g$)}\\
&\big| \sum_{i=1}^{m_{k+1}-1}g_i(x_{s_{k+1}})- \sum_{i=1}^{m_{k+1}-1}g_i(x_o)  |< \delta \qquad &\text{(by (\ref {e:<delta}))}
\end{alignat*}
Thus condition (i) is satisfied for $k+1$. Since $G_{k+1}$ is continuous,  choose $t_{k+1}>t_k$ so to satisfy (ii). By the convergence of $\sum_{i=1}^\infty g_i(x)$ for all $x$, choose $m_{k+2}>n_{k+1}$ so to satisfy (iii). 

Thus by induction we can continue the construction for all $k$ and obtain the function $G:=\lim G_k$.  As a sum of nonnegative continuous functions, $G$ is nonnegative lower semicontinuous. Then $F=g-G= \sum _{j=1}^\infty\sum_{i=n_j+1}^{m_{j+1}-1}g_i$, hence $F$ too is nonnegative lower semicontinuous.
Since
$G(x_{s_k})\ge G_k(x_{s_k}) >G_k(x_o)+\delta$ for all $k$, we have 
$$\varliminf_k G(x_{s_k})\ge \lim_k G_k(x_o) +\delta=G(x_o)+\delta.$$ Since $x_{s_k}\to x_o$, it follows that $G$ is not continuous at $x_o$.

By (iii), $$0\le G(x_{t_k})-G_k(x_{t_k})= \sum_{j=k+1}^\infty \sum_{i=m_j}^{n_j}g_i(x_{t_k})\le \sum_{i=m_{k+1}}^\infty g_i(x_{t_k})< \frac{\delta}{k}$$  and similarly $0\le G(x_o)-G_k(x_o)< \frac{\delta}{k}$.  Since by (ii),  $|G_k(x_{t_k})- G_k(x_o)|<  \frac{\delta}{k}$, it follows that 
$$|G(x_{t_k})- G(x_o)|<  \frac{3\delta}{k}$$
and hence $G(x_{t_k})\to G(x_o)$. 
Then set $y_k:= x_{t_k}$. Since $g(y_k)>\beta > g(x_o)$, it follows that $F(y_k)
\not\to F(x_o)$.   

Finally, if $K$ is convex and $g\in \LAffk$, then by \cite[Proposition 11.8]{Goodearl} 
 and  \cite[Lemma 4.2]{TikuisisToms},  $g$ is the supremum of an increasing sequence of functions in  $\Affk$ and thus we can assume that  $g_i\in \Affk_{++}$. The rest of the conclusions are now immediate.
\ep

\bL{L: chain below} 
Let $K$ be a metrizable Choquet simplex, $h\in \LAffk_{++}$, and assume there is a sequence $\{x_n\}_1^\infty \subset \Extk$ of distinct elements for which $\lim_n h(x_n)= \infty$. Then $h$ can be decomposed into the sum of two functions $F$ and $G\in \LAffk_{++}$ such that 
\begin{itemize}
\item [(i)] $G(x_n) <\infty$ for every $n$
\item [(ii)] $\sup_n G(x_n) =\infty$
\item [(iii)] $\sup_n\frac{h(x_n)}{ G(x_n)}=\infty$. 
\end{itemize}
\eL

\bp
By Proposition  \ref {P: sup aff} (i)  there is an increasing  sequence of functions $h_m\in \Affk_{++}$  that converges pointwise to $h$. Start with integers $n_1>1$ such that $h(x_{n_1})> 1$ and $m_1\ge 1$ such that $h_{m_1}(x_{n_1})\ge 1$. Then construct recursively two strictly increasing sequences of integers $n_k$ and $m_k$ such that
\be{e: ineq} h_{m_k}(x_{n_k})\ge k^2+ k\|h_{m_{k-1}}\|_\infty.\ee 

Let $\gamma:= \underset{x\in K} {\min }\,h_{m_1}(x)$. Since $ h_{m_1}$ is strictly positive, it follows that $\gamma > 0$.
Let $X_1:=\{x_1,x_2, \cdots, x_{n_1}\}$  and define for every $x_j\in X_1$
$$
g_{1,0}(x_j ):=\begin{cases}\frac{1}{2}\gamma&1\le j < n_1\\
\frac{1}{2} h_{m_1}(x_{n_1})& j= n_1.
\end{cases}
$$
We verify the conditions of Theorem \ref {T:GoodearlHB}: $X_1$ is a compact subset of $\Extk$, $g_{1,0}\in C(X_1, \mathbb R)$, the constant function $\frac{1}{2}\gamma$ is continuous and convex on $K$,  the function $\frac{1}{2}h_{m_1}\in \Affk$ is continuous and concave on $K$ and $$\frac{1}{2}\gamma\mid _{X_1}\le g_{1,0}\le \frac{1}{2} h_{m_1}\mid _{X_1}.$$ Thus by Theorem \ref {T:GoodearlHB} there is an extension $g_1\in \Affk$ of $g_{1,0}$ for which $$\frac{1}{2}\gamma\le g_1\le \frac{1}{2} h_{m_1}.$$ In particular, $g_1\in \Affk_{++}$. Let $f_1:= h_{m_1}- g_1$. Then $f_1\ge \frac{1}{2} h_{m_1}$ and hence also $f_1 \in \Affk_{++}$.

Now for every $k>1$, set $X_k:=\{x_1,x_2, \cdots, x_{n_k}\}$
 and define for every $x_j\in X_k$
$$
g_{k,0}(x_j ):=\begin{cases}0&1\le j < n_k\\
\frac{1}{k} \big(h_{m_k}-h_{m_{k-1}}\big)(x_{n_k})& j= n_k.
\end{cases}
$$
Then $g_{k,0}\in C(X_k, \mathbb R)$ and $0\le g_{k,0}\le \frac{1}{k} \big(h_{m_k}-h_{m_{k-1}}\big)\mid_{X_k}$ hence it has an extension $g_k\in \Affk$ for which \be{e: gk}0\le g_k\le \frac{1}{k}\big(h_{m_k}-h_{m_{k-1}}\big)~\text{and}~  g_k(x_{n_k})= \frac{1}{k} \big(h_{m_k}-h_{m_{k-1}}\big)(x_{n_k}).\ee Set $f_k:= h_{m_k}-h_{m_{k-1}}-g_k$. Then $f_k\ge 0$ and $f_k\in \Affk$. 

Set $m_o:=0$ and $h_{m_o}=0$. Then for all $k$
$$\sum_{j=1}^k f_j+ \sum_{j=1}^k g_j= \sum_{j=1}^k (h_{m_j}-h_{m_{j-1}})= h_{m_k}.$$
In particular,
\be{e:  Gk less}  \sum_{j=1}^k g_j< h_{m_k}\quad \forall k.\ee
Let $F=\sum_{k=1}^\infty f_k$ and $G=\sum_{k=1}^\infty g_k$, then $F+G=h$ and $F,\,G\in \LAffk_{++}$, where the strict positivity of $F$ and $G$ follows from the strict positivity of $f_1$ and $g_1$.

To show that (i) holds, for every $n$, choose $n_k>n$. By definition, $g_{k'}(x_n)=0$ for every $k'>k$ and hence by (\ref {e:  Gk less}), $$G(x_n)=\sum_{i=1}^kg_i(x_n)\le h_{m_k}(x_n)< \infty.$$
Now 
\begin{alignat*}{2}  G(x_{n_k})&\ge g_k(x_{n_k})\\
&= \frac{1}{k} \big(h_{m_k}-h_{m_{k-1}}\big)(x_{n_k})&\text{(by (\ref {e: gk}))}\\
&\ge \frac{1}{k}h_{m_k}(x_{n_k})-\frac{1}{k}\|h_{m_{k-1}}\|_\infty\\
&\ge  \frac{1}{k}(k^2+k \|h_{m_{k-1}}\|_\infty) -\frac{1}{k}\|h_{m_{k-1}}\|_\infty\qquad\qquad &\text{( by (\ref {e: ineq}))}\\
&> k
\end{alignat*}
and hence (ii) holds.
Finally 
\begin{alignat*}{2} 
G(x_{n_k})&=\sum_{j=1}^{k-1} g_j(x_{n_k})+ g_k(x_{n_k})\\
&\le h_{m_{k-1}}(x_{n_k}) + \frac{1}{k} \big(h_{m_k}-h_{m_{k-1}}\big)(x_{n_k})\quad&\text{(by (\ref {e:  Gk less}) and (\ref {e: gk}))}\\
&\le \|h_{m_{k-1}}\|_\infty + \frac{1}{k} h_{m_k}(x_{n_k})\\
&\le \frac{2}{k} h_{m_k}(x_{n_k})&\text{(by (\ref {e: ineq}))}\\
&\le \frac{2}{k} h(x_{n_k})
\end{alignat*} 
whence (iii) follows.
\ep

 \section{Quasicontinuous scale and ideals in $\Ma$}

Kucerovsky  and Perera introduced in  \cite{KucPer} the notion of quasicontinuity of the scale for simple C*-algebras of real rank zero in terms of quasitraces. In \cite {KNZCompPos} we studied this notion in terms of traces.   
\bD{D:quasicontinuous scale}\cite [Definition 2.10]{KNZCompPos}
Let $\A$ be a simple C*-algebra  with nonempty tracial simplex $\TA$.  The scale  $\Sc:=\widehat{1_{\Ma}}$  of $\A$ is said to be quasicontinuous if:
\begin{itemize}
\item [(i)]  the set $F_\infty:=\{\tau\in \Ext\mid \Sc(\tau)= \infty\}$ is finite (possibly empty) and hence the face $\co( F_\infty)$ is closed;
\item [(ii)] the complementary face $F_\infty'$ of $\co( F_\infty)$ is closed (possibly empty);
\item [(iii)] the restriction $\Sc\mid_{F_\infty'}: F_\infty'\to (0, \infty]$ of the scale $\Sc$ to  $F_\infty'$ is continuous and hence finite-valued.
\end{itemize}
\eD

As we have remarked in \cite[after Definition 2.10] {KNZCompPos}, while the scale function $\Sc$ depends on the normalization chosen  for $\TA$, the quasicontinuity of $\Sc$ does not. Notice also that when $|\Ext|< \infty$,  the scale is necessarily quasicontinuous.  If $\A$ is the stabilization of a unital algebra and hence $\Sc(\tau)=\infty$ for all $\tau\in \TA$, then $F_\infty= \Ext$ and  thus the scale is quasicontinuous if and only if $|\Ext|< \infty$. Algebras with quasicontinuous scale have interesting regularity properties. Among them, and essential for the main result of this paper is:
\bT{T:strict compar}\cite[Theorem 6.6]{KNZCompPos}
Let $\A$ be a  $\sigma$-unital simple C*-algebra with quasicontinuous scale and with strict comparison of positive elements by traces. Then strict comparison of positive element by traces holds in $\Ma$.
\eT
Extending the work by Lin \cite [Theorem 2] {LinIdeals}  on AF algebras to simple, non-unital, non-elementary  C*-algebras that are the stabilization of a unital algebra,  have strict comparison of positive elements by traces, and have a finite tracial extremal boundary, R{\o}rdam \cite [Theorem 4.4] {RordamIdeals} proved that their multiplier algebras have only finitely many ideals ($2^m-1$ when $m= |\Ext|$).  In a related  result, Kucerovsky and Perera proved (\cite [Corollary 3.5]{KucPer}) for the case of simple, separable,  non-unital, non-elementary  C*-algebras, with real rank zero, stable rank one, strict comparison of positive elements by quasitraces,  and quasicontinuous scale, that there are finitely many ideals  in $\Ma$.

The techniques in \cite{KNZCompPos} permit us to extend these results to algebras with quasicontinuous scale.  
\bT{T:f inters ideals}  Let $\A$ be a simple, $\sigma$-unital, non-unital,  C*-algebra, with quasicontinuous scale and strict comparison of positive elements by traces. For every $B\in \Ma_+\setminus \A$, let $T(B): =\{\tau\in F_\infty \mid B\in I_\tau\}$ and $I(B)$ be the ideal generated by $B$. If  $T(B)\ne \emptyset$, then $I(B)=\underset{\tau\in T(B)} {\bigcap} I_\tau$; if $T(B)=\emptyset$, then $I(B)=\Ma$.
\eT

To prove  Theorem \ref {T:f inters ideals},  we need  the following theorem and two lemmas obtained in \cite{KNZCompPos}. For the convenience of the readers and ease of reference we reproduce them here.

\bT{T:bidiag}\cite [Theorem 4.2]{KNZCompPos}
Let $\A$ be a $\sigma$-unital $C^*$-algebra and let  $T \in \Mul(\A)_+$. Then for every $\epsilon
> 0$ there exist a bi-diagonal series  $\sum_1^\infty d_k$ with each $d_k\in \A_+$ and a selfadjoint element $t_{\epsilon} \in \A $
with $\| t_{\epsilon}\| < \epsilon$ such that 
$T = \sum_1^\infty d_k + t_{\epsilon}.$  The elements $d_k$ can be chosen in $\Ped$.

For every approximate identity  $\{e_n\}$ of $\A$ with $e_{n+1}e_{n} = e_n$, we can choose $d_k$ and $t_{\epsilon}$ that satisfy the above conditions and such that for every $n\in \mathbb N$ there is an $N\in \mathbb N$ for which $e_n\sum_N^\infty d_k=0$.
\eT

 For the next lemma, notice that in \cite{KNZCompPos} we did set  $F(B)= \co\{\tau\in F_\infty\mid B\not \in I_\tau\}$ and then $T(B)= F_\infty \setminus (F(B)\cap \Ext).$

\bL{L:Lemma 5.1}\cite [Lemma 5.1] {KNZCompPos}
Let $\A$ be a  simple, $\sigma$-unital, non-unital,  C*-algebra, with strict comparison of positive elements by traces, let $a_i, \,b_i\in \A_+$  be such that $\sum_{i=1}^\infty a_i$ and $\sum_{i=1}^\infty b_i$ are two bi-diagonal series in $\Mul(\A)_+$,  let $F$ be a closed face of $\TA$, $F'$ be its complementary face  (either $F$ or $F'$  can be empty), and assume that $|F\cap\Ext|< \infty$.
Assume  also that  
for some $\epsilon, \delta, \alpha >0$ we have  
\item [(i)]$ \big(\sum\nolimits_{i=1}^\infty b_i-\delta\big)_+\not \in \A$
\item [(ii)] $d_{\tau}\Big(\big(\sum\nolimits_{i=m}^\infty b_i-\delta\big)_+\Big)=\infty \hspace{5.3cm}\forall \, \tau \in F, \, m\in\mathbb N$,
\item [(iii)] $d_{\tau}\Big(\big(\sum\nolimits_{i=1}^\infty a_i-\epsilon\big)_+\Big)+\alpha \le d_{\tau}\Big(\big(\sum\nolimits_{i=1}^\infty b_i-\delta\big)_+\Big)< \infty \hspace{2.05cm}\forall~\tau\in F'$, 
\item [(iv)] $d_{\tau}\Big(\big(\sum\nolimits_{i=m}^n b_i-\delta\big)_+\Big)\to d_{\tau}\Big(\big(\sum\nolimits_{i=m}^\infty b_i-\delta\big)_+\Big)$ \hspace{0.4cm} uniformly on $F'$,  $\forall\,m\in \mathbb N$,
\item [(v)] $d_{\tau}\Big(\big(\sum\nolimits_{i=n}^\infty a_i-\epsilon\big)_+\Big)\to 0$ \hspace{5.25cm}uniformly on $F'$.

Then  
$\Big(\sum_{i=1}^\infty a_i-2\epsilon\Big)_+\preceq\Big (\sum_{i=1}^\infty b_i-\delta'\Big)_+$ for all $\delta'$ with $0< \delta'< \delta.$
\eL
\bL{L:Lemma 6.4}\cite [Lemma 6.4]{KNZCompPos}
 Let $\A$ be a simple, $\sigma$-unital, non-unital,  C*-algebra,  $P \in \Mul(\A)$ be a projection,
$K\subset \TA$ be a closed set such that
$\widehat{P}\mid_K$ is continuous, and let
$ \sum_{j=1}^{\infty} A_j$ be the strictly converging sum of elements $A_j\in (P \Mul(\A) P)_+$. Assume furthermore that there exists an increasing approximate identity 
$\{ e_n \}_{n=1}^{\infty}$ for $(P \A P)_+$
with $e_{n+1} e_n = e_n$ for all $n \in \mathbb N$
such that for all $m \geq 1$, there exists  $N  \in \mathbb N$ with
$e_m \sum_{j = N}^{\infty} A_j = 0$.
Then for every $\delta \ge 0$, 
\item [(i)]
$d_{\tau}\Big( \Big( \sum_{j=n}^{\infty} A_j - \delta 
\Big)_+ \Big) \to 0$ uniformly on $K$.
\item [(ii)] $d_{\tau}\Big( \Big( \sum_{j=1}^n A_j - \delta \Big)_+ \Big) \to  d_{\tau}\Big( \Big(  \sum_{j=1}^{\infty} A_j - \delta \Big)_+ \Big)$ uniformly on $K$.
\eL
The above two lemmas are based on the following result which we also will need in our paper:

\bP{P:diag Cuntz} \cite[Proposition 4.4]{KNZCompPos}
Let $\A$ be a C*-algebra, $A=\sum_1^\infty A_n$, $B=\sum_1^\infty B_n$ where $A_n, B_n\in \Mul(\A)_+$, $A_nA_m=0$, $B_nB_m=0$ for $n\ne m$ and the two series converge in the strict topology. If $A_n\preceq (B_n-\delta)_+$ for some $\delta >0$ and for all $n$, then $A\preceq (B-\delta')_+$ for all $0< \delta'< \delta$. 
\eP
The proof of Theorem \ref {T:f inters ideals}, which is based on the above two lemmas, is inspired by the proof of  \cite[Theorems 5.3 and 6.6]{KNZCompPos}.

\bp[Proof of Theorem \ref {T:f inters ideals}]
Assume that $T:=T(B)\ne \emptyset$, i.e., $B\in I_\tau$ for some  $\tau\in F_\infty$, leaving to the reader the similar (and simpler) case when $T=\emptyset$. Set $$F=\co\{F_\infty\setminus T\}.$$ Then \be{e:F finite} |F\cap\Ext|\le  |F_\infty|< \infty. \ee Being finite dimensional, the face  $F$ is closed and hence split, i.e., $\TA=F\overset{\cdot}{+}F' $ where $F'$ is  the complementary face of $F$. Then
\be{e:F'} F'= \co(T)\overset{\cdot}{+}F'_\infty,\ee  
is also closed since $F'_\infty$ is closed by hypothesis and  $ \co(T)$ is closed  because it is finite dimensional. 
 Since $I(B)\subset  \bigcap\{I_\tau\mid \tau\in T\}$, we need to prove that if $A\in \Ma_+$ and $A\in I_\tau$ for all $\tau\in T$, then $A\in I(B)$.
We can assume that $\|A\|=\|B\|=1$ and by using \ref {T:bidiag}, we reduce to the case that $A=\sum_{k=1}^\infty a_k$ and $B=\sum_{k=1}^\infty b_k$ are bidiagonal series for an approximate identity $\{e_n\}_{n=1}^\infty$ and  for all $m \geq 1$, there exists  $N  \in \mathbb N$ with
$e_m \sum_{k = N}^{\infty} a_k = 0$ and $e_m \sum_{k = N}^{\infty} b_k = 0$. Since $A$ decomposes into the sum of two diagonal series  $A= \sum_{k=1}^\infty a_{2k-1}+\sum_{k=1}^\infty a_{2k}$, to simplify notations we can assume that $A$ is diagonal. Choose $\delta>0$ such that 
\be{e:delta}
(B-\delta)_+\not \in \A\quad\text{and}\quad (B-2\delta)_+\not \in I_\tau\quad \forall \tau\in F_\infty\setminus T. \ee
 Let $\epsilon >0$.
Since 
$$(A-\frac{\eps}{2})_+ +\big(I_{\Ma}- (A-\frac{\eps}{2})_+\big)=I_{\Ma},$$ $\widehat{(A-\frac{\eps}{2})_+}$ is complemented under the scale $\Sc$ and hence it is continuous on $F_\infty'$. As it is continuous also on the finite dimensional face $ \co(T)$, it follows that 
$$ \widehat{(A-\frac{\eps}{2})_+}= \sum_{k=1}^\infty  \widehat{(a_k-\frac{\eps}{2})_+}\in \mathrm{Aff}(F')_{++}.$$
By Dini's theorem, the series  $\sum_{k=1}^\infty  \widehat{(a_k-\frac{\eps}{2})_+}$ converges uniformly on $F'$. Let $$2\alpha:=\min \{d_\tau\big((B-\delta)_+\big)\mid \tau\in \TA\}.$$ Choose $N$ such that $\sum_{k=N}^\infty  \widehat{(a_k-\frac{\eps}{2})_+}(\tau)\le \frac{\eps}{2}\alpha$ for all $ \tau\in F'.$
Then
$$d_\tau\big(\sum_{k=N}^\infty a_k-\eps\big)_+= \sum_{k=N}^\infty d_\tau( a_k-\eps)_+\le \frac{2}{\eps}\sum_{k=N}^\infty \tau( a_k-\frac{\eps}{2})_+\le \alpha \quad \forall \,\tau\in F'$$
and thus
\be{e:cond iii}
d_\tau\big(\sum_{k=N}^\infty a_k-\eps\big)_++\alpha\le 2\alpha\le d_\tau\big((B-\delta)_+\big)\quad \forall \tau\in F'.
\ee

Now we are in the position to verify that all the hypotheses (i)-(v) of Lemma \ref {L:Lemma 5.1} are satisfied for the diagonal series $A_N= \sum_{k=N}^\infty a_k$, the bidiagonal series $B= \sum_{k=1}^\infty b_k$, the face $F$, and the scalars $\eps, \delta$, and $\alpha$. 

By  (\ref {e:delta}),  the hypothesis (i) of Lemma  \ref {L:Lemma 5.1} holds  and  also  
$(B-2\delta)_+\not\in I_\tau$ for every $\tau\in F$. Since by (\ref {e:ineq 3})
$$
d_\tau\Big(\big(\sum_{k=1}^\infty b_k -2\delta\big)_+\Big)\le d_\tau\Big(\big(\sum_{k=1}^{m-1} b_k -\delta\big)_+\Big)+ d_\tau\Big(\big(\sum_{k=m}^\infty b_k -\delta\big)_+\Big)
$$
and since $$d_\tau\Big(\big(\sum_{k=1}^{m-1} b_k -\delta\big)_+\Big)\le \frac{2}{\delta}\tau\Big(\big(\sum_{k=1}^{m-1} b_k -\frac{\delta}{2}\big)_+\Big)< \infty$$ as $\sum_{k=1}^{m-1} b_k \in \A$, it follows that
\be{e:cond (ii)}
d_\tau\Big(\big(\sum_{k=m}^\infty b_k -\delta\big)_+\Big)= \infty\quad \forall \tau\in F, m\in \mathbb N\ee
which establishes hypothesis (ii) of Lemma  \ref {L:Lemma 5.1}. Hypothesis (iii) was  established in (\ref {e:cond iii}). 
Since $\widehat{1_{\Ma}}=\Sc$ is by hypothesis continuous on $F_\infty'$ and since 
$\sum_{k=N}^\infty a_k$ converge strictly, it follows from Lemma \ref {L:Lemma 6.4} that $d_\tau\Big(\big (\sum_{k=m}^\infty a_k-\epsilon\big)_+\Big) \to 0$ uniformly on $F'_\infty$. By \cite [Lemma 3.2] {KNZCompPos}, 
$d_\tau\Big(\big (\sum_{k=m}^\infty a_k-\epsilon\big)_+\Big) \to 0$ for every trace $\tau\in T$.  Since $ T$ is finite,  $d_\tau\Big(\big (\sum_{k=n}^\infty a_k-\epsilon\big)_+\Big) \to 0$ uniformly  also on $\co(T)$ and hence by  (\ref {e:F'}) the convergence is uniform also on $F'$. 

By the same argument, for every $m\in \mathbb N$, the strict convergence of 
$$\big (\sum_{k=m}^n b_k-\delta\big)_+ \to \big (\sum_{k=m}^\infty b_k-\delta\big)_+$$ implies the uniform convergence over $F'$ of 
$$d_\tau\Big(\big (\sum_{k=m}^n b_k-\delta\big)_+\Big) \to d_\tau\Big(\big (\sum_{k=m}^\infty b_k-\delta\big)_+\Big).$$
Thus  conditions (v) and (iv) of Lemma \ref  {L:Lemma 5.1} are also established. Therefore,
$$
\Big(\sum_{k=N}^\infty a_k-2\epsilon\Big)_+\preceq \sum_{k=1}^\infty b_k $$
and hence $\Big(\sum_{k=N}^\infty a_k-2\epsilon\Big)_+\in I(B)$. Since 
$$(A-2\eps)_+= \Big(\sum_{k=1}^\infty a_k-2\epsilon\Big)_+= \Big(\sum_{k=1}^{N-1}a_k-2\epsilon\Big)_++\Big(\sum_{k=N}^\infty a_k-2\epsilon\Big)_+$$ and 
$\Big(\sum_{k=1}^{N-1}a_k-2\epsilon\Big)_+\in \A\subset I(B)$, it follows that $(A-2\eps)_+\in I(B)$. As $\eps$ is arbitrary, we conclude that $A\in I(B)$.
\ep

As a consequence we obtain:
\bC{C: structure ideals}
 Let $\A$ be a simple, $\sigma$-unital, non-unital,  non-elementary, C*-algebra, with strict comparison of positive elements by traces, and with quasicontinuous scale, and let $m:=|F_\infty|$. Then $\Ma$ has precisely $2^m-1$ proper ideals properly containing $\A$, each being an intersection of  ideals  $I_\tau$ for $\tau\in F_\infty$.
\eC
Notice that if $\A=\K$ then $m=1$ but there are no proper ideals 
properly containing  $\A$, thus for the exact count of the ideals in $\Ma$ we need indeed to assume that $\A$ is non-elementary. 

\section{Projection-surjectivity and injectivity}\label{S: proj inj surj}
We find it convenient to introduce the following  terminology for properties that have appeared in various forms in the study of multiplier algebras of C*-algebras.
\bD{D:surject and inject} Let $\A$ be a simple, $\sigma$-unital, non-unital,  C*-algebra, with non empty tracial simplex $\TA$.
\item [(i)] $\A$ is 1-projection-surjective if for every $f\in \LAffs_{++}$ that is complemented under $\Sc=\widehat{1_{\Ma}}$ (i.e., there is $g\in  \LAffs_{++}\sqcup\{0\}$ such that $f+g=\Sc$) there is a projection $P\in \Ma\setminus \A$ such that $f=\hat P$. 
\item [(ii)] $\A$ is 1-projection-injective if  $P\sim Q$  whenever $P, \,Q\in \Ma\setminus \A$ are projections such that $\hat P= \hat Q$.
\item [(iii)] $\A$ is n-projection-surjective (resp., n-projection-injective) if the algebra $M_n(\A)$ is 1-projection-surjective (resp., 1-projection-injective).
\item [(iv)] $\A$ is projection-surjective and injective if it is 1-projection-surjective and 2-projection-injective.
\eD

Notice that $\K$ is obviously not 1-projection-surjective, thus whenever we assume 1-projection-surjectivity it is redundant to require that the algebra be non-elementary.
We start with some simple relations between n-projection-surjectivity and  m-projection-injectivity for various m and n.

\bL{L:prj inj.surj relations}
 Let $\A$ be simple, $\sigma$-unital,  non-unital, C*-algebra, with nonempty and metrizable tracial simplex $\TA$.
 \item [(i)]  If $\A$ is n-projection-surjective for some $n\in \mathbb N$, then it is kn-projection-surjective for every  $k\in \mathbb N$ and $\St$ is 1-projection surjective. 
  \item  [(ii)] If $\A$ is n-projection-injective (resp., $\St$ is 1-projection-injective),  
 then it is k-projection-injective for every $k< n$ (resp., every $k\in \mathbb N$). 
 \item [(iii)]  Let $\A$ be n-projection-surjective (resp., $\St$ is 1-projection-surjective). If $\A$ is 2n-projection-injective (resp.. $\St$ is 1-projection-injective), then $\A$ is 1-projection-surjective.
 \item [(iv)] If $\St$ is 1-projection-injective and surjective, then $\A$ is n-projection-injective and surjective for every $n$.

\eL 
\bp
\item [(i)] Assume that $\A$ is n-projection-surjective, let $k\in \mathbb N$,  and let $f+g= kn \Sc$ for some $f\in \LAff_{++}$ and $g \in \LAff_{++}\sqcup\{0\}$. Then $\frac{f}{k}$ is complemented under $n\Sc$ and hence there is a projection $P\in \Mul(M_n(\A))\setminus M_n(\A)= M_n(\Ma)\setminus M_n(\A)$ such that $\hat P=\frac{f}{k}$. Then $Q:=\bigoplus_{j=1}^k P\in M_{kn}(\Ma)\setminus M_{kn}(\A)$ and $\hat Q = k\hat P= f$.  Thus $\A$ is kn-projection-surjective.
We prove now that $\St$ is 1-projection-surjective. Let  $f\in \LAff_{++}$. By the metrizability of $\TA$ and Proposition \ref {P: sup aff},  $f=\sum_{j=1}^\infty f_j$ with $f_j\in \Aff_{++}$. For every $j$, choose $n_j> \frac{\max f_j}{\min \Sc}$ and $n_j$ divisible by $n$. Then $f_j< n_j \Sc$ and since $f_j$ is continuous, $f_j$ is complemented under $n_j \Sc$. By the first part of the proof, $\A$ is $n_j$-projection-surjective, hence there is a projection $P_j\in \Mul(M_{n_j}(\A))\setminus M_{n_j}(\A)$ such that $\widehat {P_j}= f_j$. Construct a strictly converging series of mutually orthogonal projections $\tilde{P_j}$ in $\M$ such that $\tilde{P_j}\sim P_j$ and the series $P=\sum_{j=1}^\infty \tilde{P_j}$ converges strictly. Then $P\not \in \St$ and $$\hat P=\sum_{j=1}^\infty \widehat{\tilde{P_j} }= \sum_{j=1}^\infty \widehat{P_j}= \sum_{j=1}^\infty f_j = f.$$
\item [(ii)] Assume that $\A$ is n-projection-injective, let $k\le n$ and let $P,Q$ be projections in $\Mul(M_k(\A))\setminus M_k(\A)$ with $\hat P=\hat Q$. Then $P\oplus 0$, $Q\oplus 0$ are projections belonging to  $\Mul(M_n(\A))\setminus M_n(\A)$ and $\widehat {P\oplus 0}= \widehat {Q\oplus 0}$. Then $P\oplus 0\sim Q\oplus 0$ and hence $P\sim Q$.
\item [(iii)]  Assume that $\A$ is n-projection-surjective and 2n-projection-injective and let $f \in \LAff_{++}$  and $g\in \LAff_{++}\sqcup\{0\}$ such that $f+g= \Sc$. Then $f+ (n-1)f+ng=n\Sc$, i.e., $f$ is complemented under $n\Sc$ and so is $g$. Thus there are projections $P,Q\in M_n(\Ma)$ with $\hat P=f$ and $\hat Q= g$. Hence there are  mutually orthogonal  projection $P', Q'\in M_{2n}(\Ma)$ with $P'\sim P$, $Q'\sim Q$. But then $$\widehat {P'\oplus Q'} =\hat P+\hat Q= f+g= \Sc= \widehat{1_{\Ma}}.$$ Since $M_{2n}(\A)$ is  1-projection-injective by hypothesis, $P'+Q'\sim 1_{\Ma}$. Thus we can choose $P', Q'$ with $P'+Q'= 1_{\Ma}$ and hence $P', Q'\in \Ma$. In particular, $\hat{P'}=f$. \\
The case when  $\St$ is 1-projection-surjective and 1-projection-injective is similar and is left to the reader.
\item [(iv)] Obvious.
\ep

In all cases where we could determine projection-surjectivity and injectivity, the property holds for every $n$. Does 1-projection-injectivity imply 2-projection-injectivity and hence n-projection-injectivity for every $n$?  The answer is affirmative in the case when the algebra has real rank zero.
\bL{L:RR0 2x2}
Let $\A$ be a $\sigma$-unital, non-unital, non-elementary, simple C*-algebra of real rank zero and let $P\in M_2(\Ma)$ be a projection. Then $P\sim P''\oplus P''$ for some projection $P''\in \Ma$.
\eL
\bp
Let $\{e_n\}$ be an increasing approximate identity for $\A$ consisting of projections and set $e_0=0$. Then $\{e_n\oplus e_n\}$ is an increasing approximate identity of projections  for $M_2(\A)$. By \cite [Theorem 4.1]{SZ90-2}, and passing if necessary to a subsequence of $\{e_n\}$, we can find projections $p_n\le (e_n-e_{n-1})\oplus (e_n-e_{n-1}) $ such that $P\sim  \sum _{n=1}^\infty p_n $ and  the series converges in the strict operator topology. By \cite [Theorem 3.3]{SZ90-2},  we can further assume that for all $n$,
$$p_n=\begin{pmatrix}s_n&0\\0&s_n+r_n\end{pmatrix}$$ for some projections in $\A$,  $s_n, r_n\le e_n- e_{n-1}$. By a slight adjustment of the proof, we  can  also  assume that  $r_n \ne 0$ for all $n$. 
By \cite [Theorem 1.1] {ZhangMatricial} we can approximately halve $r_1$, that is decompose it into the sum of three mutually orthogonal projections
$r_1= t_1+t_1'+  q'_2$
where $t_1\sim t_1' $ and $q'_2\sim q_2\lneq  r_2$. 
Then since  $$ \begin{pmatrix}0&0\\0&t'_1\end{pmatrix}\sim \begin{pmatrix}t_1&0\\0&0\end{pmatrix}$$ and both  are orthogonal to $\begin{pmatrix}s_1&0\\0&s_1+ t_1\end{pmatrix}$, it follows that
$$\begin{pmatrix}s_1&0\\0&s_1+ t_1+ t_1'\end{pmatrix}\sim \begin{pmatrix}s_1+t_1&0\\0&s_1+ t_1\end{pmatrix}.$$
Similarly, 
$$\begin{pmatrix}0&0\\0&q'_2\end{pmatrix}\sim  \begin{pmatrix}0&0\\0&q_2\end{pmatrix}\sim\begin{pmatrix}q_2&0\\0&0\end{pmatrix}$$
and 
$$
p_1=  \begin{pmatrix}s_1&0\\0&s_1+ t_1+ t_1'+ q'_2\end{pmatrix}\sim \tilde p_1:=\begin{pmatrix}s_1+t_1&0\\0&s_1+t_1\end{pmatrix}+ \begin{pmatrix}q_2&0\\0&0\end{pmatrix}.
$$

Next, approximately halve $r_2-q_2= t_2+t_2'+  q'_3$ with $ t_2\sim t_2'$ and $q'_3\sim q_3\lneq r_3$. Then reasoning as above,
$$p_2=  \begin{pmatrix}s_2&0\\0&s_2+ t_2+ q_2+ t_2'+ q'_3\end{pmatrix}\sim p'_2:= \begin{pmatrix}s_2+ t_2&0\\0&s_2+ t_2+q_2\end{pmatrix}+ \begin{pmatrix}q_3&0\\0&0\end{pmatrix}.$$
Since $p_1p_2=0$ and $p_1'p_2'=0$, it follows that $p_1+p_2\sim p_1'+p_2'$, namely
$$
p_1+p_2\sim \begin{pmatrix}s_1+t_1&0\\0&s_1+t_1\end{pmatrix}+\begin{pmatrix}s_2+t_2+q_2&0\\0&s_2+t_2+q_2\end{pmatrix}+ \begin{pmatrix}q_3&0\\0&0\end{pmatrix}.
$$

Iterating, we find a sequence of  mutually orthogonal projections $$s_n,t_n, q_n\le e_n-e_{n-1}$$ such that for every $n$
$$p_n\sim p_n':=\begin{pmatrix}s_n+ t_n&0\\0&s_n+ t_n+q_n\end{pmatrix}+ \begin{pmatrix}q_{n+1}&0\\0&0\end{pmatrix}.$$
 Since $p_n'\le e_{n+1}-e_{n-1}$, the series  $P':=\sum_{n=1}^\infty p'_n$ converges strictly. Choose partial isometries  $v_n\in M_2(\A)$ such that $p_n= v^*_nv_n$ and $p'_n= v_nv^*_n$. Then the series $V:=\sum_{n=1}^\infty p'_nv_np_n$ also converges strictly to the partial isometry $V\in M_2(\Ma)$. Then $P=V^*V$,  $P'= VV^*$, and hence $P\sim P'$ within $M_2(\Ma).$ 
 Setting $q_1:=0$ we have for every $k$ that
$$\sum_{n=1}^k p'_n = \sum_{n=1}^k\begin{pmatrix}s_n+t_n+ q_n&0\\0&s_n+t_n+ q_n\end{pmatrix} + \begin{pmatrix}q_{k+1}&0\\0&0\end{pmatrix}$$
and hence
$$
P' = \sum_{n=1}^\infty \begin{pmatrix}s_n+t_n+ q_n&0\\0&s_n+t_n+ q_n\end{pmatrix}.$$
Let $P''=  \sum_{n=1}^\infty s_n+ t_n+ q_n$, then $P''\in \Ma$ is a projection and   $P'= P''\oplus P''$, which completes the proof.

\ep 
\bP{P:RRO 1inj=2inj}Every separable, non-unital, non-elementary, simple C*-algebra of real rank zero which is 1-projection-injective is also n-projection-injective for all $n\ge 1$.
\eP
 \bp
 In view of Lemma \ref {L:prj inj.surj relations}, it is sufficient to prove the statement for $n=2$. Assume that $P, Q$ are projections in $ \Mul(M_2(\A))$ and that $\hat P=\hat Q$. By Lemma \ref {L:RR0 2x2}, $P\sim P''\oplus P''$, $Q\sim Q''\oplus Q''$ for some projections $P''$ and $Q''$ in $\Ma$.  Hence $\hat {P''}=\hat {Q''}$ and hence $P''\sim Q''$ whence  $P\sim Q$.
 
 \ep

We proceed now to ascertain projection-surjectivity and injectivity for some important classes of simple, $\sigma$-unital, non-unital,  non-elementary, C*-algebras.
 
We start with the case of real rank zero algebras with stable rank one which was long well-known (\cite{ElliottHandelman}, \cite {Zhang89}, \cite{LinNotes}, \cite{Lin1}, \cite {GoodearlK0}, \cite {LinQuasidiagonal}). A nice exposition can be found in \cite [Theorem 3.9] {PereraIdeals}. 
\bT{T:Perera}  Let $\A$ be a simple, $\sigma$-unital, non-unital,  non-elementary C*-algebra, with real rank zero, stable rank one, and such that $\A$ has strict comparison of positive element by traces. Then $\A$ is n-projection-surjective and n-projection-injective for every $n$. 
\eT
\bp The hypotheses in  \cite [Theorem 3.9] {PereraIdeals} on the C*-algebra $\A$ are that $\A$ is simple, $\sigma$-unital, non-unital,  non-elementary, has  real rank zero, stable rank one, and that  the monoid $V(\A)$ of equivalent classes of projections in $M_\infty(\A)$ is strictly unperforated. The latter hypothesis is equivalent to the condition that $\A$ has strict comparison of positive elements by 2-quasitraces (see  Lemma 3.5, Corollary 3.10 and its proof in \cite {PereraStructure}). Obviously, strict comparison of positive elements by traces implies strict comparison by quasitraces, so the hypotheses of  \cite [Theorem 3.9] {PereraIdeals} are satisfied. The thesis of  \cite [Theorem 3.9] {PereraIdeals} is expressed in terms of a monoid isomorphism of $V(\Ma)$ which implies $n$-projection surjectivity and injectivity of $\A$ for every $n$. 

\ep 

The condition that $\A$ has real rank zero can be dropped in the case when $\A$  is separable and stable. 
\bT{T:Huaxin and Ping 1} \cite[Proposition 4.2]{LinNg} Let $\A$ be a simple, non-unital, separable, C*-algebra, with stable rank one, and with strict comparison of positive elements by traces. Then $\A$ is 1-projection-injective.
\eT
\cite[Proposition 4.2]{LinNg} requires the algebra to be stable, but an examination of its proof shows that stability is not necessary.
Next we consider projection-surjectivity.
\bT{T:Huaxin and Ping 3}  \cite[Corollary 4.6]{LinNg} Let $\A$ be a simple, separable C*-algebra with non empty tracial simplex $\TA$ such that \\
(*) for every bounded function $f\in \LAff_{++}$ there exists an  $a\in (\St)_+$  which is not Cuntz equivalent to a projection and such that $\widehat{[a]} = f$. \\
Then $\St$ is 1-projection-surjective. 
\eT
The property (*) in the above theorem plays an important role in the study of Cuntz semigroups and is succinctly  formulated in \cite {BPT08} as the surjectivity of the map $\iota:W(\A)_+\mapsto \LAffb_{++}$  where $W(\A)_+$ is the sub-semigroup of equivalence classes of elements of $M_\infty(\A)_+$ not equivalent to projections and $\iota[a](\tau)= d_\tau(a)$ for all $a\in  M_\infty(\A)_+$ and $\tau\in \TA$. The property (*) was first shown to hold for C*-algebras that are simple, unital, separable and are either exact, stably finite, and $\mathscr Z$-stable (\cite[ Theorem 5.5.]{BPT08}) or are infinite-dimensional AH algebras of stable rank one  with strict comparison of positive elements (\cite[ Theorem 5.3.]{BPT08}). Condition (*) is also satisfied by some stably projectionless algebras, e.g., the monotracial Razak algebra (see for example \cite {RobCone}.)

The $\mathscr Z$-stability condition in (\cite[ Theorem 5.5.]{BPT08}) was recently replaced by the weaker condition of having stable rank one. 
 
\bT{Thiel Th. 8.11}\cite [Theorem 8.11]{ThielRanks} Let $\A$ be a separable, unital, simple, non-elementary C*-algebra with stable rank one. Then for every $f\in \LAffq_{++}$ there exists $x\in \St_+$ such that $d_\tau(x)=f(\tau)$ for all $\tau\in \mathcal QT(\A)$.
\eT

Here $\mathcal QT(\A)$ denotes the Choquet simplex of 2-quasitraces of $\A$, which contains $\TA$ as a face. Notice that the hypothesis can be reformulated by asking $\A$ to be stable and to contain a non-zero projection.

The statement of this theorem does not state explicitly that $x$ can be chosen to be not equivalent to a projection. It is easy to see that this can be done when $\A$ has real rank zero:

\bP{P: RRO} 
Let $\A$ be a  simple, real rank zero C*-algebra with non-empty  tracial simplex $\TA$. Then for every projection $q\in \A_+$ there exists $a\in \A_+$ such that $a\le q$, $a$ is not equivalent to a projection, and $\widehat{[a]}=\hat q$.
\eP

\bp
We can assume that $\A$ is nonelementary, as the elementary case is trivial. By \cite [Theorem 1.1] {ZhangMatricial}, we can decompose $q$  into a sum of projections $q= q_1+q_1'$ with $0\ne q_1'\preceq q_1$ and hence with $\widehat{q_1'}\le \frac{1}{2}\hat q$. Then decompose similarly $q_1'= q_2+q_2'$ with $0\ne \widehat{q_2'}\le \frac{1}{2^2}\hat q$. Iterating the process, we find an infinite sequence of mutually orthogonal projections $q_n\le q$ such that 
$$q- \sum_{n=1}^m q_n= q_m'\quad\text{and}\quad  \widehat{q_m'}\le \frac{1}{2^m}\hat q.$$ 
Since $\hat q$ is continuous, it follows that $\hat q= \sum_{n=1}^\infty \hat q_n$. Then $a:=\sum_{n=1}^\infty \frac{1}{n}q_n\in \A_+$, $a\le q$, $\widehat {[a]}= \hat q$ and $a $ is not equivalent to a projection because $0$ is an accumulation point in the spectrum of $a$.
\ep

The same result holds also for (stable, separable) algebras that don't have real rank zero due to the work  \cite{APT2018},  presented in \cite [Proposition 2.9] {ThielRanks}, that states that for a countably based, simple, stably finite, non-elementary Cuntz semigroup $S$ satisfying axioms (05) and (06) (and hence for the concrete Cuntz semigroup of the stable C*-algebra $\A$ considered) for every $[a]\in S$ there is $[x]\in S$, $[x]\le [a]$ and $[x]$ {\it soft} (and hence $x$ is not equivalent to a projection) such that $\tau (x)=\tau(a)$ holds for all 2-quasitraces and hence a fortiori for all traces $\tau$. We summarize this result for our setting:

\bP{P: Proposition 2.9} . Let $\A$ be a separable, simple, non-elementary, stable C*-algebra with non-empty  tracial simplex $\TA$. Then for every $x\in \A_+$ there exists $a\in \A_+$ such that $a\preceq x$, $a$ is not equivalent to a projection, and $\widehat{[a]}= \widehat{[x]}$
\eP

Combining Theorem \ref {Thiel Th. 8.11}, Proposition \ref {P: Proposition 2.9}, and    Theorem \ref {T:Huaxin and Ping 3} we obtain:
\bC{C:surjectivity}
Let $\A$ be a separable, unital, simple, non-elementary C*-algebra with stable rank one. Then $\St$ is 1-projection surjective. If furthermore $\St$ has strict comparison of positive elements by traces, then $\St$ is also 1-projection-injective.
\eC 

Thus the class of C*-algebras $\A$ with both projection injectivity and projection surjectivity  for $\St$ includes among others :
}
\begin{itemize}
\item real rank zero algebras with stable rank one and strict comparison of positive elements, including all simple unital AF-algebras and all irrational rotation algebras,
\item all simple finite nuclear C*-algebras that have been classified in the Elliott program,
\item all crossed products of the form 
$C(X) \times_{\alpha} \mathbb{Z}$, where
$X$ is a compact metric space with finite topological dimension and
$\alpha : X \rightarrow X$ is a minimal homeomorphism, 
\item the Jiang--Su algebra $\mathcal Z$ and more generally, all  simple, unital, separable, exact, stably finite $\mathcal Z$-stable C*-algebras,
\item The reduced free group C*-algebra $C_r^*(F_\infty)$  on infinitely many generators, 
\item The monotracial Razak algebra (stably projectionless).
\ \end{itemize}

Notice that all the C*-algebras listed above also have strict comparison of positive elements (by traces). We will prove that under the additional hypothesis of separability and stability, strict comparison is indeed necessary for projection-surjectivity and injectivity. We need first a simple consequence of the definition of projection-surjectivity and injectivity and of the argument in the proof of Lemma \ref {L:prj inj.surj relations} (iii) that will be useful throughout the rest of the paper. 
\bL{L:subequiv} Let $\A$  be  a simple, $\sigma$-unital, non-unital, non-elementary,  C*-algebra, and let $P$, $Q$ be projections in $\Ma$. 
\item [(i)] If  $P\preceq Q$ then $\hat P$ is complemented under $\hat Q$. \\
Assume now  that $\A$ is projection-surjective and injective 
and that $Q\not\in \A$. 
\item [(ii)] If $f+g=\hat Q$ for some $f, g\in \LAffs_{++}$, then there is a decomposition of $Q=  P_1+P_2$  into projections $P_1, P_2\not\in \A$ with $\hat{ P_1}= f$ and $\hat{ P_2}= g$. 
\item [(iii)] If  $\hat P$ is complemented under $\hat Q$, then $P\preceq Q$.  
\eL
\bp
\item [(i)] There is a projection $P'\in \Ma$ with  $P\sim P'\le Q$ and hence $\hat P=\hat {P'}$. Let $P''= Q-P'$, then $\hat Q= \hat {P}+\hat{P''}$ and since $\hat{P''}$ is either 0 (if $P'' =0$) or strictly positive (if $P''\ne 0$),  it follows that $\hat P$ is complemented under $\hat Q$.
 \item [(ii)] Since 
 $f+ g+ \widehat {(1_{\Ma}- Q)}= \widehat {1_{\Ma}}=\Sc$
 both $f$  and $g$ are complemented under $\Sc$. Thus there are projections $R_1, R_2\not\in\A$ such that $\hat  {R_1}=f$ and $\hat  {R_2}=g$. Then $\widehat {R_1\oplus R_2}= g+ f= \hat Q$. Since neither $R_1\oplus R_2\in M_2(\A)$ nor $Q\oplus 0\in M_2(\A)$, by 2-projection injectivity, $R_1\oplus R_2\sim Q\oplus 0$ and hence $Q=P_1+P_2$ for some mutually orthogonal projections  $P_1\sim R_1$ and $P_2\sim R_2$. Thus $\hat  {P_1}=f$ and $\hat  {P_2}=g$.
 \item [(iii)] Let $g\in \LAffs_{++}$ be such that $\hat P+g= \hat Q$. Reasoning as in the proof of (ii), there is a projection $R_2\in \Ma\setminus \A$ such that $\widehat {P\oplus R_2}=  \hat Q$. Since neither $P\oplus R_2$ not $Q$  are in $\A$, it follows that $P\oplus R_2\sim Q\oplus 0$ and hence $P\prec Q$.
\ep


Next, we list the following facts that are routine, but by completeness we add a short proof.
\bL{L: subequiv vs strict}
Let $\B$ be a C*-algebra.  
\begin{enumerate} \item [(i)] Let $T\in \Mb_+$, $T_n\in \Mb_+$ such that $T_n\to T$ strictly. If $a\in \B_+$ and $a\preceq T$, then for every $\eps>0$  there is an $n$ such that  $(a-\eps)_+\preceq T_n$.
\item  [(ii)]
Let $Q\in \Mb$ be a projection and assume that $Q\B Q$ has a strictly positive element $b$. If $a\in \B_+$ and $a\preceq Q$, then $a\preceq b$.
\end{enumerate}
\eL
\bp
(i)
Choose an $X\in \Mb$ such that $\|a- XTX^*\|< \eps/3$, an  $e\in \B_+$ with $\|e\|=1$ such that $\|a-eae\|< \eps/3$, and  an integer $n$ such that $\|eX(T-T_n)X^*e\|< \eps/3$. Then
$\|a- eXT_nX^*e\|< \eps$ and hence $$(a-\eps)_+\preceq eXT_nX^*e\preceq T_n.$$
(ii) Let $\eps>0$. Since $b^{1/n}$ converges strictly to $Q$, by (i) there is an integer $n$ such that $(a-\eps)_+\preceq b^{1/n}\sim b$. Since $\eps $ is arbitrary,  then $a\preceq b$.

\ep

We need also a standard application of Kasparov's Absorption Theorem which has appeared in many places over the years (e.g., \cite{KucNg}, \cite {LinNg}). The precise form of the argument that we require can be found in Lemma 4.3 and the proof of Proposition 4.4 in \cite {LinNg}.

\bL{L: Kasparov}
$\A$ be a simple, stable, separable C*-algebra and let $a\in \A_+$. Then there is $a'\in \A_+$ with $a\sim a'$ and $R_{a'}\in \Ma$. Furthermore, $R_{a'}\in \A$ if and only if $a$ is equivalent to a projection.
\eL
\bp
By Kasparov's Absorption Theorem and \cite [Lemma 4.3]{LinNg}, there is a projection $P\in \Ma$ such that the Hilbert modules $\overline {a\A}$ and $P\A$ are isomorphic, i.e., there there is a unitary $\Phi: \overline {a\A}\mapsto P\A$. If $b$ is a strictly positive element in $\A$, then $a':=PbP$ is a strictly positive element in $P\A P$ and $R_{a'}= P$. Moreover, $P\A= \overline {a'\A}$. Then by a standard argument (see for example \cite[Proposition 4.3] {OrtegaRordamThiel}, see also \cite {CowardElliottIvanesu}, \cite {LinCuntz}), $a\sim a'$. If $P\in \A$, then $a'\sim P$ and hence $a$ is equivalent to a projection. Conversely, if $a$ equivalent to a projection $P$, then we can choose $a'=P$.  \ep

\bT{T: proj surj inj implies strict comp}
Let $\A$ be a simple, stable, separable,  C*-algebra with projection-surjectivity and injectivity. Then $\A$ has strict comparison of positive elements by traces.
\eT
Notice that by our definition, projection-surjectivity or injectivity implies that $\A$ has non-empty tracial simplex and, clearly, projection-surjectivity implies that $\A$ is non-elementary.
\bp
Let $a, b\in \A_+$ and assume that $d_\tau(a)< d_\tau(b)$ for every $\tau\in \TA$ such that $d_\tau(b)<\infty$.  Assume without loss of generality that $\|a\|\le 1$. By Lemma \ref {L: Kasparov} we can assume that $R_a\in \Ma$. 
For the first step of the proof, we construct for every $\eps>0$ a projection $P\in \Icon\setminus \A$, such that  $(a-\eps)_+\preceq P$ and $\hat P< \widehat {[b]}$.

Consider first the case when $R_a\not \in \A$, namely when $a$ is not equivalent to a projection. Since $\widehat{R_a}\in \LAffs_{++}$, by Proposition \ref {P: sup aff} (i) we can decompose $\widehat{R_a}$ into the pointwise converging sum $\widehat{R_a}= \sum_{n=1}^\infty f_n$ of functions $f_n\in \Aff_{++}$. 
By projection-surjectivity, we can find projection $R_n''\in \Ma\setminus\A$ such that $\widehat {R_n''}=f_n$ for every $n$. Since $\A$ is stable, we can find mutually orthogonal projections $ R_n'\sim R_n''$ such that $R'=\sum _{n=1}^\infty R_n'$ converges strictly. 
As $\widehat R'= \sum_{n=1}^\infty f_n= \widehat{R_a}$, by projection-injectivity we have $R_a\sim R'$. This provides a strictly converging decomposition of $R_a=\sum_{n=1}^\infty R_n$ into projections $R_n \in \Icon\setminus \A$.  
Let $\eps>0$. Then by Lemma \ref  {L: subequiv vs strict}, there is an $n$ such that
 $$(a-\eps)_+\preceq P:= \sum_{k=1}^nR_k.$$Thus $P \in\Icon \setminus \A$,  and $\hat P< \widehat {R_a}=  \widehat{[a]}\le  \widehat {[b]}$.  
 
Next consider the case when $R_a\in \A$. Then  $ \widehat{[a]}= \widehat {R_a}$ is continuous, hence  $$\widehat{[b]}-\widehat {R_a}\in \LAff_{++}\quad\text{and}\quad \Sc- \widehat {R_a}\in \LAff_{++}.$$  Let $0< \alpha< \min \big(\widehat {[b]}-\widehat {R_a}\big)$. 
Then also $\Sc- \widehat {R_a}-\alpha\in  \LAff_{++}$, hence the constant function $\alpha$ is complemented under $
\Sc- \widehat {R_a}=\widehat{1_{\Ma}-R_a}$. By Lemma \ref {L:subequiv}, there is a projection $P_o\in \Ma\setminus \A$ with $P_o\le I_{\Ma}-R_a$ and $\widehat{P_o}= \alpha$. Let $P=R_a+P_o$. Then $P \in\Icon \setminus \A$, 
$(a-\eps)_+\le a\le R_a\le P$ and $\hat P< \widehat{[b]}.$

For the second step of the proof, by \cite [Proposition 2.9] {ThielRanks}, there is $c'\in \A_+$, $c'\preceq b$, with $\widehat{[c']}= \widehat{[b]}$ and $c'$ not equivalent to a projection. Again by Lemma \ref {L: Kasparov}, there is a $c\in \A_+$ with $c\sim c'$ and such that $R_c\in \Ma\setminus \A$. Since $\hat P< \widehat {R_c}= \widehat {[b]}$ and $\hat P$ is continuous, it follows that $P\preceq R_c$ (see Corollary  \ref {C: strict comp in Icon}  below). As $c$ is strictly positive in $R_c\A R_c$, it follows by Lemma \ref {L: subequiv vs strict} that $(a-2\eps)_+\preceq c\preceq b$. As $\eps$ is arbitrary, it follows that $a\preceq b$.

\ep
If $\A$ is just $\sigma$-unital and/or if $\A$ is not stable, we cannot invoke Proposition \ref {P: Proposition 2.9}.  However, if  $\A$ has real rank zero, we can still prove strict comparison for $\A$.

\bP{P: strict comp RR0}
Let $\A$ be simple, $\sigma$-unital, non-unital C*-algebra with real rank zero and with projection-surjectivity and injectivity. Then $\A$ has strict comparison of positive elements by traces.
\eP

\bp It is well-known (e.g., see \cite [Corollary 3.10] {PereraStructure} and its proof) that it suffices to prove that A has strict comparison of projections by traces. Let $p$, $q$ be projections in $\A$, and assume that $\hat p(\tau) <\hat q(\tau)$ for all $\tau\in \TA$. Since $\hat p$  and $\hat q$ are continuous, $\hat q- \hat p$ is continuous and $\Sc- \hat p$ is lower semicontinuous.  Choose $$0< \alpha < \min\{\hat q(\tau)- \hat p(\tau)\mid \tau\in \TA\}.$$  Then the constant function $\alpha$ is complemented under $\Sc- \hat p= \widehat{ 1_{\Ma}-p}$. By Lemma \ref {L:subequiv}, there is a  $P_o\in \Ma\setminus \A$ orthogonal to $p$ and such that $\widehat{P_o}=\alpha$. Thus $P:=p+P_o\in \Ma\setminus \A$, $\hat P$ is continuous, and $\hat P(\tau) < \min \{q(\tau)\mid \tau\in \TA\}$.

Reasoning as in the proof of Proposition \ref {P: RRO}   we can find a sequence of mutually orthogonal nonzero projections $q_n\le q$ such that $\hat q = \sum_{n=1}^\infty \widehat{q_n}$. By Dini's theorem the convergence is uniform, so there is $N$ such that $\sum_{n=1}^N\widehat{q_n}> \hat P$. To simplify notations, assume that $N=1$. Now choose an approximate identity $ \{e_n\}$ of $\A$ consisting of projections  and such that $e_1=q_1$. Since $\A$ is simple and of real rank zero, we can find for every $n\ge 2$ projections  $0\ne q'_n\le e_n-e_{n-1}$ and $q_n'\sim q_n''\le q_n$ 
 Set $q'_1:=q_1$ and $Q':=\sum_{n=1}^\infty q'_n$. Since the series converges strictly, $Q' \in \Ma\setminus \A$ and  $$\widehat {Q'}> \widehat{q_1}> \hat P.$$
$\hat P$ being continuous, it is complemented under $\widehat Q'$. By Lemma \ref {L:subequiv}, $P\preceq Q'$ and hence $p\preceq Q'$. By Lemma \ref {L: subequiv vs strict}, there is a $n$ such that
$$p\sim (p-\frac{1}{2})_+\preceq \sum_{k=1}^n q'_k\sim \sum_{k=1}^n q''_k\le \sum_{k=1}^n q_k\le q,$$
which completes the proof.
\ep

\section{Projection-surjectivity and injectivity and ideals in $\Ma$}\label{S: first appl}
As this section will illustrate, assuming that a C*-algebra  is projection-surjective and injective greatly facilitates the study of the ideal structure of its multiplier algebra. 

\bP{P:dividing} Let $\A$  be  a simple, $\sigma$-unital, non-unital, C*-algebra,  which is projection-surjective and injective and let  $P\in \Ma\setminus \A$ be a projection. Then for every $n\in \mathbb N$ there are mutually orthogonal projections $P_1\sim P_2\cdots \sim P_n$ in $\Ma$ such that $P=\sum_{j=1}^nP_j$.\eP 
\bp
Since $\hat P= \frac{1}{n}\hat P+  \frac{n-1}{n}\hat P$ and both functions are in $\LAffs_{++}$, by Lemma \ref {L:subequiv} there are  a mutually orthogonal projection $P_1$ and $P_1'$ not in $\A$
such that $Q= P_1+ P_1'$,  $\widehat {P_1}= \frac{1}{n}\hat P$, and $\widehat {P_1'}= \frac{n-1}{n}\hat P$. By the same reasoning,   $P_1'$ is  the sum of two orthogonal projections $P_1'= P_2+P_2'$ not in $\A$ such that  $P_2\sim P_1$ and  $\widehat{P_2'}= \frac{n-2}{n}\hat P$. After $n-1$ steps we get  a decomposition of $P$ into mutually orthogonal projections not in $A$, $P= P_1+\cdots +P_{n-1}+ P_{n}$ with $P_1\sim P_2\sim \cdots \sim P_{n-1}$  and with $\widehat{P_{n}}= \frac{1}{n}\hat P$. Then $P_n\sim P_1$ by projection-injectivity,  which completes the proof. 
\ep
Compare this result with the case when $\A$ has real rank zero  where it was shown in \cite {ZhangMatricial} that projections in $\Ma\setminus \A$ are divisible  by $2^m$.  Notice that as a consequence, for every $n$, $\Ma\simeq M_n(\Mb)$ for some hereditary subalgebra $\B\subset \A$.

\bC{C:V(M(A))}
Let $\A$  be  a simple, $\sigma$-unital, non-unital, C*-algebra,  which is n-projection-surjective and injective for all $n$ and let  $[P]\in V(\Ma)\setminus V(\A)$. Then there is an $n\in \mathbb N$ such that $[P]= \sum_{j=1}^n[P_j]$ for some projections $P_j\in \Ma\setminus \A$. \eC
\bp
Since $P\in M_n(\Ma)= \Mul(M_n(\A))$ for some $n\in \mathbb N$, and $P\not \in V(\A)$ and hence in particular, $P\not \in M_n(\A)$, by Proposition \ref {P:dividing}, $P= \sum_{j=1}^nP_j$ with $P_j\sim P_1\in M_n(\Ma)$. Then $\hat P $ is complemented under $n\Sc$, i.e., there is a function $f\in \LAffs_{++}\sqcup \{0\}$ such that
$\hat P+f=n\Sc$. But then $\widehat{P_1}+\frac{f}{n}= \Sc$ and hence by Lemma \ref {L:subequiv}, $\widehat{P_1}=\widehat{P_1'}$ for some projection $P_1'\in \Ma\setminus \A$. By n-projection-injectivity, $P_1\sim P_1'$ and hence the conclusion follows. 
\ep

Another  simple consequence of Lemma \ref {L:subequiv} is 
\bP{P: in ideal}  Let $\A$  be  a simple, $\sigma$-unital, non-unital,   C*-algebra, which is projection-surjective and injective, and let $P$ and $Q$ be projections in $\Ma$ with  $Q\not \in \A$. Then the following conditions are equivalent.
\item [(i)] $P\in I(Q)$ (the principal ideal generated by $Q$);
\item [(ii)]  $\hat P + f=m \hat Q$ for some  $m\in \mathbb N$ and some $f\in \LAffs_{++}\sqcup \{0\}$.
\eP
\bp
Assume that (i) holds, then for some $m\in \mathbb N$, by Lemma \ref {L:ideals}, $$P\preceq \bigoplus_{k=1}^m Q\in M_m(\Ma)= \Mul(M_m(\A)).$$ Hence by Lemma \ref {L:subequiv}(i), there is an $f\in \LAffs_{++}\sqcup \{0\}$ such that $$\hat P+ f= \widehat {\bigoplus_{k=1}^m Q}=m\hat Q.$$\\
Assume that (ii) holds, then by Proposition \ref {P:dividing} we can decompose 
$P$ into the sum of $m$ mutually orthogonal and equivalent projections, $P=\sum _{k=1}^m P_k$ and hence $\widehat {P_k}= \frac{1}{m}\hat P$ for every $k$. Then
$\widehat {P_k}+\frac{1}{m}f = \hat Q$. By Lemma \ref {L:subequiv} (ii), $P_k\preceq Q$ and  hence $P_k\in I(Q)$ for every $k$, whence  $P\in I(Q)$.
\ep

In the case when $\A$ is simple, $\sigma$-unital, non-unital, non-elementary, and has strict comparison of positive elements by traces, we proved in  \cite [Theorem 6.4]{KNZMin} that strict comparison of positive elements holds for $\Icon$.   The proof depended on the technique developed in \cite{KNZCompPos} and used in the present paper in Theorem \ref {T:f inters ideals}. As the following corollary illustrates, in the presence of projection-surjectivity and injectivity, strict comparison of projections for $\Icon$ can be obtained with a considerably simpler proof and without requiring explicitly strict comparison for the underlying algebra $\A$ (which however holds automatically by Theorem \ref {T: proj surj inj implies strict comp} if we further assume that $\A$ is separable and stable).
\bC{C: strict comp in Icon}   Let $\A$  be  a simple, $\sigma$-unital, non-unital,   C*-algebra,  which is projection-surjective and injective, and let $P\in \Icon$ and  $Q\in \Ma\setminus  \A$ be projections. 
\item [(i)]  If $\hat P(\tau)<\hat Q(\tau)$ for all $\tau$, then $P\preceq Q.$
\item [(ii)] If $P\not\in\A$ and $Q\in \Icon$, then $I(P)=I(Q)$.
\eC
\bp
\item [(i)] Since $\hat P \in \Aff_{++}$ by (\ref {e:cont proj}), it follows that $f:= \hat Q- \hat P\in \LAff_{++},$ and hence  $\hat P$ is complemented (by $f$) under $ \hat Q$.  Thus   $P\preceq Q$ by Lemma  \ref {L:subequiv} (iii). 
\item [(ii)] Since also $\hat Q \in \Aff_{++}$  we can 
choose $n$ such that $ \max \hat P< n \min \hat Q$.  By Proposition \ref {P:dividing} ,  decompose $P$ into the sum of $n$ mutually orthogonal equivalent projections, $P=\sum_{k=1}^n P_k$, with $\widehat {P_k} =\frac{1}{n}\hat P$. Every $\widehat {P_k}$ is continuous, hence $P_k\in \Icon$. By (i), $P_k\preceq Q$ for every $k$, hence $P_k \in I(Q)$  and thus $P\in I(Q)$. Interchanging the role of $P$ and $Q$, we conclude that  $I(P)=I(Q)$. 

\ep

In  \cite [Theorem 6.6]{KNZCompPos} we proved that if $\A$ is $\sigma$-unital, simple, has quasicontinuous scale, and has strict comparison of positive elements, then $\Ma$ has strict comparison of positive elements (see Definition \ref {D:str comp}) (see also \cite {KNZCompJOT} for the real rank zero case).  In the presence of projection-surjectivity and injectivity, Corollary \ref {C: strict comp in M(A)} here below will show that strict comparison of projections for $\Ma$ can be obtained much more easily and without requiring explicitly strict comparison for the underlying algebra $\A$.   We will use the notation introduced in Theorem \ref {T:f inters ideals} for a projection $P\in \Ma$:
$$T(P) = \{\tau\in F_\infty \mid \tau(P)< \infty\}.$$

\bC{C: strict comp in M(A)} Let $\A$  be  a simple, $\sigma$-unital, non-unital,   C*-algebra, which is projection-surjective and injective and has  quasicontinuous scale, and  let $P, Q\in \Ma$ be projections with $Q\not \in \A$.
\item [(i)] If  $\hat P(\tau)<\hat Q(\tau)$ for all $\tau$ such that $\hat Q(\tau)< \infty$, then  $P\preceq Q.$ 
\item [(ii)] If $P\not \in \A$ and $T(P)=T(Q)$, then $I(P)=I(Q)$.
\eC

\bp
\item[(i)] Set  $T:= T(Q)$  and $F:= \co(F_\infty\setminus T).$
We will assume that  $T\ne\emptyset$ and $F\ne  \emptyset$, as the case when one of the two sets is empty is similar but simpler and will be left to the reader.
The face $F$ is finite dimensional, and thus  closed (and hence split) (\ref{e: closed face  is split}). Its complementary face  $F'$ itself splits as
$F'= \co(T) \overset{\cdot}{+} F_\infty'$. As $F'$ is the direct sum of the  finite dimensional and hence closed face  $\co(T)$ and the face $F_\infty'$ which is closed by hypothesis, it is also closed.
Since $\hat Q+ \widehat {I-Q}= \Sc$, $\hat P+ \widehat {I-P}= \Sc$ and $ \Sc $ is continuous on $F_\infty'$, by Lemma \ref{L:facts}, $\hat Q$ and $\hat P$ are continuous on $F_\infty'$. Both functions are also continuous on $\co(T)$ since $\hat P(\tau)<\hat Q(\tau)< \infty$ for every $\tau\in   T$ and $T$ is finite. Thus  $\hat Q-\hat P\in \mathrm{Aff}( F')_{++}$. 
Then by Lemma \ref {L:direct sum} (iii), $$f:= \hat Q\mid_{F}\overset{\cdot}{+} (\hat Q- \hat P)\mid_{F'}\in \LAff_{++}.$$ Since $$\hat P(\tau)+ f(\tau)=\begin{cases} \hat P(\tau)+\hat Q(\tau)=\infty &\tau\in F\\\hat Q(\tau)& \tau\in F'\end{cases}\quad=\quad\hat Q(\tau),$$
$\hat P$ is complemented under $\hat Q$ and hence $P\preceq Q$ by Lemma \ref {L:subequiv}.
\item [(ii)] By the first part of the proof, both $\hat Q$ and $\hat P$ are continuous on the closed face $F'$. Thus we can find $n$ such that $$\underset{\tau \in F'}{\max}(\frac{1}{n}\hat P(\tau))<\underset{\tau \in F'}{ \min }\hat Q(\tau).$$
Reasoning as in the proof of Corollary \ref {C: strict comp in Icon} (ii), $P=\sum_{k=1}^n P_k$ and $\widehat {P_k}(\tau)< \hat Q(\tau)$ for $\tau\in F'$, i.e., for all $\tau$ such that $\hat Q(\tau)< \infty$. By part (i), $P_k\preceq Q$ and hence $P\in I(Q)$. Interchanging the role of $P$ and $Q$, we conclude that  $I(P)=I(Q)$. 
\ep

\bC{C:order ideals}
Let $\A$ be a simple, $\sigma$-unital, non-unital,  C*-algebra with quasicontinuous scale, and metrizable tracial simplex $\TA$, and which is  n-projection-surjective and n-projection injective for every integer $n$. Let $H$ be an order ideal of $V(\Ma))$ not contained in $V(\A)$. Let $$S:=\{\tau\in F_\infty\mid \exists [P]\in H\text{ such that } \tau(P)=\infty\}.$$ Then 
$H= \{ [Q]\in V(\Ma))\mid \tau(Q)< \infty \text{ for all } \tau \in F_\infty\setminus S\}$.
In particular, $V(\Ma))$ has only finitely many order ideals.
\eC

\bp
By definition, $H\subset \{ [Q]\in V(\Ma))\mid \tau(Q)< \infty \text{ for all } \tau \in F_\infty\setminus S\}$. To prove the opposite inclusion, for every $\tau\in S$, choose  $[P_\tau]\in H$ such that $\tau(P_\tau)=\infty$. Let 
$P:= \underset{\tau\in S}\bigoplus P_\tau$. 
Since $S\subset F_\infty$ is finite it follows   that   $[P]\in H$ and that $\tau(P)= \infty$ for all $\tau\in S$. Since $P \in M_k(\Ma)$ for some $k$, it is complemented under $k\Sc$ and reasoning as in the proof Corollary \ref {C: strict comp in M(A)},  it is continuous on $\co(S)'= \co(F_\infty\setminus S)\+F_\infty'$. 
Let $[Q]\in V(\Ma))$ be such $ \tau(Q)< \infty$  for all $\tau \in F_\infty\setminus S$.  By the same reasoning as for $\hat P$, $ \hat Q$ is continuous on $\co(S)'$.  By Lemma \ref {C:complementation} (iii), $\hat Q $ is complemented under $m\hat 
 P$ for some integer $m$. By the assumption of n-projection-surjectivity and  injectivity for every $n$ and by Lemma \ref {L:subequiv}, it follows that $[Q]\le m[P]$ and hence $[Q]\in H$.
 
 \ep

As a further consequence of projection surjectivity and injectivity we obtain the maximality for the ideals $I_\tau$ when $\tau\in F_\infty=\{\tau\in \Ext\mid \Sc(\tau)= \infty\}.$ Maximality for $I_\tau$ was obtained for the stable case by R{\o}rdam \cite [Theorem 4.4]{RordamIdeals} for $\A\otimes \K$, $\A$ unital, with strict comparison of positive elements by traces 
and finite extremal boundary. The same result was also obtained by Perera  in the proof of \cite [Theorem 6.6]{PereraIdeals} for quasitraces  and $\sigma$-unital, non-unital, non-elementary C*-algebras with real rank zero, stable rank 1, and weakly unperforated $K_0$ group.  These results generalized earlier work by \cite {ElliottMatroidII}, \cite {LinMultIdeals}.
\bT{T:maximal in stable}
Let $\A$ be a simple, separable, non-unital,   C*-algebra,  such that $\St$ is projection-surjective  and injective  and let $\tau_o\in F_\infty$. 
\item [(i)] The ideal $I_{\tau_o}$ of $\Ma$ is generated by any projection $P\not \in \A$ such that \\$\hat P(\tau):\begin{cases}< \infty& \tau=\tau_o\\
=\frac{\Sc}{2} & \tau\in\{\tau_o\}'\end{cases}$. Such projections exist.
\item [(ii)] $I_{\tau_o}$ is a maximal ideal.
\eT
\bp
We first prove both these statements under the additional hypothesis that $\A$ is stable, in which case $\Sc(\tau)=\infty$ for all $\tau\in \TA$ and $F_\infty=\Ext$. Notice that by separability of $\A$,  $ \LAffs_{++}=  \LAff_{++}$ by Proposition \ref {P: sup aff} (i) and every function in $ \LAff_{++}$ is complemented under $\Sc$. 
\item [(i)] 
Let $g: =1\mid_{\{\tau_o\}}\overset{\cdot}{+}\frac{\Sc}{2}\mid_{\{\tau_o\}'}$, or,  more explicitly, $g(\tau)=\begin{cases}1& \tau=\tau_o\\
\infty & \tau\ne\tau_o\end{cases}$. By Corollary \ref {C:complementation} (or directly from the definition), $g\in \LAff_{++}$. Thus by Lemma \ref {L:subequiv},   there is a projection $P\in \Ma\setminus \A$ such that $\hat P= g$.  By (\ref {e:proj in Itau}), $P\in I_{\tau_o}$.
We claim that every positive $A\in  I_{\tau_o}  $ belongs to $I(P)$. By expressing $A$ as an $\A$ perturbation of the sum of two positive diagonal elements (Theorem \ref{T:bidiag})  and then reasoning as in the proof of Theorem \ref {T:f inters ideals},   we can assume that $A$ itself is diagonal, i.e.,  $A=\sum_1^\infty a_n$, where $a_n\in \A_+$,  $a_na_m=0$ for $n\ne m$, and the series converges in the strict topology.  Fix $\eps>0$, then  by (\ref {e:ineq 2}) and (\ref {e:proj in Itau}),
$$\sum_{n=1}^\infty d_{\tau_o}((a_n-\eps)_+)=d_{\tau_o}((A-\eps)_+)\le \frac{2}{\eps}\tau_o((A-\frac{\eps}{2})_+)< \infty~ \forall ~\eps>0.$$ 
Choose $N$ such that $\sum_{n=N}^\infty d_{\tau_o}((a_n-\eps)_+)< 1$. As it is enough to prove that $A_N:= \sum_{n=N}^\infty a_n\in I(P)$, to simplify notations assume that $N=1$. 
By the stability of $\A$, decompose $I_{\Ma}= \sum_{k=1}^\infty E_k$ into a sum of mutually orthogonal projections $E_k\sim I_{\Ma}$.
Let $$\alpha_n:= d_{\tau_o}\big((a_n-\eps)_+\big) + \frac {1-d_{\tau_o}\big((A-\eps)_+\big)}{2^n},$$ so that $\sum_{n=1}^\infty \alpha_n=1$.
Then, again by Corollary \ref {C:complementation}, $$g_n:= \alpha _n\mid_{\{\tau_o\}}\overset{\cdot}{+}\frac{\Sc}{2}\mid_{\{\tau_o\}'}\in \LAff_{++}$$ and is complemented under $\Sc$. By the 1-projection-surjectivity of $\A$ there is a projection $P_n\not \in \A$ with $\widehat {P_n}= g_n$  and since $E_n\sim I_{\Ma}$, we can take $P_n\le E_n$. 
Then the series $\sum_{n=1}^\infty P_n$ converges strictly to a projection $R$ and $$\hat R= \sum_{n=1}^\infty g_n= g= \hat P.$$ Then $R\sim P$ by the 1-projection injectivity of $\St$, so assume without loss of generality that $R=P$.  By the separability of $\A$,  for every $n$ we can find a strictly positive element $b_n\in P_n\A P_n$. Then $d_\tau(b_n)= \widehat {P_n}(\tau)$ for all $\tau$. Since
$$  d_{\tau}( a_n-\eps)_+ < \begin {cases}
\alpha_n= g_n(\tau_o)=\widehat{P_n}(\tau_o)&\tau= \tau_o\\
\infty =g_n(\tau) = \widehat{P_n}(\tau)&\tau\ne \tau_o
\end{cases} \qquad= d_\tau(b_n), $$
by Theorem \ref {T: proj surj inj implies strict comp}, we obtain from the strict comparison for $\A$ that
$( a_n-\eps)_+\preceq b_n$. 

Now $b_n\le \|b_n\|P_n\sim P_n\sim (P_n- \frac{1}{2})_+ $  and hence $ ( a_n-\eps)_+\preceq (P_n- \frac{1}{2})_+$ for every $n$. By Proposition \ref {P:diag Cuntz}, $(A-\eps)_+\preceq P$.  Since $\eps>0$ is arbitrary, $A\preceq P$ and hence $A\in I(P)$. This proves that $I_{\tau_o}= I(P).$

\item [(ii)]
Let $P$ be a projection for which $I_{\tau_o}=I(P)$, $\mathcal J$ be a closed two-sided ideal of $\Ma$ such that $I_{\tau_o}\subsetneq \mathcal J$, and let $A\in \mathcal J_+\setminus  I_{\tau_o}$.

Invoking  Theorem \ref {T:bidiag} and reasoning as in the first  part of the proof, we can assume that $A=\sum_{n=1}^\infty a _n$ with $a_n\in \A_+$ mutually orthogonal and $A\not \in I_{\tau_o}$. Choose $\eps>0$ such that $(A-\eps)_+\not\in  I_{\tau_o}$. As a consequence $$d_{\tau_o}\big((A-\eps)_+\big)= \sum_{n=1}^\infty d_{\tau_o}\big((a_n-\eps)_+\big)= \infty.$$ Let $e_n$ be an approximate identity of $\A$ such that $e_{n+1}e_n=e_n$ for all $n$ and all $\tau\in \TA$. Recall that all $e_n$ are in the Pedersen ideal of $\A$, and by Lemma \ref {L:bounded dim}, $d_\tau(e_n)< \infty$ for all $n$. By regrouping if necessary finite sums  of $a_n$ terms, assume that $d_{\tau_o}\big((a_n-\eps)_+\big) > d_{\tau_o}(e_n-e_{n-1}) $ for all $n$, where we set $e_0=0$. Reasoning as in part (i), decompose $P$ into a sum of $P=\sum_{n=1}^\infty P_n$ with 
$\begin{cases} \widehat {P_n}(\tau))< \infty &\tau=\tau_o\\\widehat {P_n}(\tau))=\infty &\tau\ne \tau_o\end{cases}$. Since $\A$ is separable, there is a strictly positive $b_n\in P_n\A P_n$ and we can assume that $\|b_n\|=1$.
Then for every $n$ and every $\tau$
$$
d_\tau(e_n- e_{n-1}) < d_\tau\big(b_n \oplus (a_n-\eps)_+\big).
$$
 Indeed, for $\tau= \tau_o$, $d_\tau(e_n- e_{n-1}) < d_\tau( (a_n-\eps)_+)$, while 
 $d_\tau(b_n )= \tau(P_n)= \infty$ for every other $\tau$.
Since $M_2(\A)$ has strict comparison, it follows that for every $n$
$$e_n- e_{n-1}\preceq b_n\oplus (a_n-\eps)_+\le P_n\oplus (a_n-\eps)_+\sim \big((P_n\oplus a_n)- \eps\big)_+.$$
Since $1_{\Ma}= \sum_{n=1}^\infty e_n- e_{n-1}$ and  $\sum_{n=1}^\infty P_n
\oplus a_n= P\oplus A$ where both series converge strictly, again by Proposition \ref {P:diag Cuntz} we obtain that $1_{\Ma} \preceq P\oplus A$. As $P\in I_{\tau_o}\subset \mathcal J$ and $A\in\mathcal J$, we have $P\oplus A\in \mathcal J$, thus $1_{\Ma}\in    \mathcal J$ and hence $\mathcal J=\Ma$.  We thus conclude that $I_{\tau_o}$ is maximal.

Finally,  we remove the hypothesis that $\A$ is stable. There is a projection $R\in \M$ such that $\A$ is isomorphic to $R(\St) R$ and hence, by identifying $1_{\Ma}$ with $R$,
$\Ma$ can be identified with $R \M R$. As usual, we identify the tracial simplex $\TA$ of $\A$ with the tracial simplex of $\St$. 
Every ideal $\mathcal J$ of $R\M R$ is the compression $\mathcal J=R\tilde{ \mathcal J}R$ of an ideal $\tilde{ \mathcal J}$ of $\M$.
For every $\tau\in \TA$, denote by $I_{\tau, \A}$ (resp., $I_{\tau, \St}$) the ideal of $R\M R$ (resp., of $\M$). It is then immediate to verify that $I_{\tau, \A}= RI_{\tau, \St}R$. Similarly, if  $P\in \M$ is a projection and $P\le R$, then $I_{\A}(P)=R I_{\St}(P)R$ where we denote by $I_{\A}(P)$ (resp., by   $I_{\St}(P)$) the principal  ideal of $R \M R$ (resp., of $\M$) generated by $P$.  Since $\tau_o\in F_\infty$, the function $$g: =1\mid_{\{\tau_o\}}\overset{\cdot}{+}\frac{\Sc}{2}\mid_{\{\tau_o\}'}\in \LAff_{++}$$ constructed at the beginning of the proof is complemented under $\Sc$ by Corollary \ref {C:complementation} and hence there is a projection $P\in R\M R$  with $\hat P=g$.  Since $\St$ satisfies the hypotheses, by the first part of the proof, $I_{\tau_o, \St}= I_{\St}(P)$ and then
$$I_{\tau_o, \A}= RI_{\tau_o, \St}R= RI_{\St}(P)R= I_\A(P).
$$
Furthermore,  since $I_{\tau_o, \St}$ is  maximal and $I_{\tau_o, \A}$ is proper, it follows that $I_{\tau_o, \A}$ is also maximal, which proves (i) and (ii) also for the case when $\A$ is not stable.

 \ep

\section{Characterization of purely infinite corona algebras}

In this section we examine the link between pure infiniteness of the corona algebra $\Ma/\A$  and other properties of the algebra $\A$ and its multiplier algebra $\Ma$. Not all the implications require the same hypotheses on the algebra $\A$. 

\bP{P:strict comp implies prop inf}
Let $\A$ be a simple, $\sigma$-unital, non-unital, non-elementary,  C*-algebra,  with non-empty tracial simplex $\TA$ and with strict comparison of positive elements by traces for $\Ma$. 
Then $\Ma/\A$ is purely infinite. 
\eP
\bp
Since no non-zero quotient of $\Ma$ can be abelian,  the corona algebra $\Ma/\A$ has no characters, hence by  \cite [Definition 4.1]{KirchRor},  to obtain that $\Ma/\A$ is purely infinite it is (necessary and) sufficient to prove that if $A, B\in \Ma_+$ and $\pi(A)\in I(\pi(B))$, then $\pi(A)\preceq \pi(B)$. Clearly, $A\in I(B)$.
By Theorem \ref{T:bidiag}, $A= \sum_1^\infty a_k +b_o$ where $b_o=b_o^*\in \A$ , $0\ne a_n\in \A_+$ and the series is bidiagonal ($a_na_m=0$ for $|n-m|>1$ and converges strictly.) Now  $a_1\sum_3^\infty a_k=0$ and $\pi(A)= \pi\big(\sum_3^\infty a_k\big)$, so to simplify notation simply assume that there is an $0\ne a\in \A_+$ such that $aA=0$. Choose a strictly positive element $b\in \A$, then for all $\tau\in \TA$ $d_\tau(b)=d_\tau(1_{\Ma}) =\Sc(\tau)$ and
$$d_\tau(A)\le d_\tau(A+a)\le d_\tau(b) = d_\tau(B+b)$$
where the first inequality is strict for all $\tau$ for which $d_\tau(b)< \infty$ and thus
$d_\tau(A)< \infty$.
Since $A\in I(B)= I(B+b)$, by the assumption of strict comparison on $\Ma$, we have $A\preceq B+b$ and hence $\pi(A) \preceq B$. 

\ep 

\bP{P: corona pi}
Let $\A$ be a simple, $\sigma$-unital, non-unital, non-elementary,  C*-algebra,  with non-empty tracial simplex $\TA$. 
Assume there exists a projection $P$ in $\Ifin$ but not in $\Icon$. Then $\pi(P)\in \Ma/\A$ is not properly infinite. In particular,  $ \Ma/\A$ is not purely infinite.
\eP
\bp
Assume by contradiction that $\pi(P)\oplus \pi(P)\preceq \pi(P)$. Then there is some $X\in  \Mul(M_2(\A))$ such that
$$ \|\pi(X)\pi(P)\pi(X)^*- \pi(P)\oplus \pi(P)\|< \frac{1}{2}$$ and hence there is some $a=a^*\in M_2(\A)$  for which
$\|XPX^*+a - P\oplus P\|< \frac{1}{2}.$ 
Let $a=a_+-a_-$ with $a_-, a_+ \ge 0$, then $(P\oplus P+a_- -\frac{1}{2})_+\preceq XPX^*+a_+$. Hence
$$ P\oplus P\sim (P\oplus P-\frac{1}{2})_+\preceq XPX^*+a_+\preceq P\oplus a_+.$$

It is well known
that then there is a $\delta>0$ and a projection $$Q\in \Her( (P\oplus a_+-\delta)_+)= \Her( P\oplus (a_+-\delta)_+)$$ such that $P\oplus P\sim Q$. Notice that $(a_+-\delta)_+$ belongs to the Pedersen ideal of $M_2(\A)$ and has also a (positive) local unit $b$ in the same Pedersen ideal, that is $b(a_+-\delta)_+=(a_+-\delta)_+$. Then $P\oplus b$ is a local unit  for $P\oplus (a_+-\delta)_+$ and hence also for $Q$, that is $(P\oplus b)Q=Q$. Thus $Q\le P\oplus b$. Let $g:= \widehat{P\oplus b-Q}$. Then $g\in \LAff_{+}$ and 
$$2\hat P+ g=\hat Q+g=  \hat P+\hat b.$$
Since $P\in \Ifin$, $\hat P(\tau)$ is finite for every $\tau \in \Ext$ and hence
$$\hat b(\tau)= \hat P(\tau)+ g(\tau)\quad\forall \tau\in \Ext.$$
By Proposition \ref {P: sup aff}, $\hat b= \hat P+ g.$
Since $\hat b$ is continuous because $b$ belongs to the Pedersen ideal and  since both functions $\hat P$ and $g$ are lower semicontinuous, it follows by Lemma \ref {L:facts} that  $\hat P$  must  be continuous. By (\ref {e:cont proj}) this contradict the hypothesis that $P\not \in \Icon$.

\ep

If $ \Ma/\A$ is purely infinite, it thus follows that all the projections of $\Ifin$ are in $\Icon$. If $\A$ is 1-projection-surjective this is sufficient to guarantee that  $\Ifin=\Icon$:  
\bL{L: Ifin not Ib}
Let $\A$ be a simple, $\sigma$-unital, non-unital C*-algebra, and assume that $\A$ is 1-projection-surjective.  
If $\Ifin\ne \Ib$ (resp.,  $\Ib\ne \Icon$), then there is a projection $P\in \Ifin\setminus \Ib$ (resp., $P\in \Ib\setminus \Icon$).
\eL
\bp
Let $A\in (\Ifin)_+\setminus \Ib$.  Without loss of generality, assume that $\|A\|\le 1$. By Lemma  \ref {L: char Ifin, Ib} there is some $\delta>0$ and some $\mu\in \TA$ for which $\widehat {(A-\delta)_+}(\mu)= \infty$. $\widehat{(A-\delta)_+}\in \LAffs_{++}$ and since $(A-\delta)_+\le I$, the evaluation function $f:= \widehat{I-(A-\delta)_+}$ also belongs to $\LAffs_{++}$. As  
$\Sc = \widehat{1_{\Ma}}= \widehat{(A-\delta)_+}+ f$, by Lemma \ref {L:subequiv}  there is a projection $P$  such that $\hat P= \widehat{(A-\delta)_+}$. As $\hat P(\tau)< \infty$ for all $\tau\in \Ext$ and $\hat P(\mu)= \infty$, it follows that $P\in \Ifin\setminus \Ib$ by Lemma  \ref {L: char Ifin, Ib}. 

The case when $\Ib\ne \Icon$ is similar: there is $A\in \Ib\setminus\Icon$ with $\|A\|\le 1$ and $\delta>0$  and a projection $P$ such that $\hat P=  \widehat{(A-\delta)_+}$ is bounded but not continuous, and hence $P\in  \Ib\setminus \Icon$ by (\ref {e:cont proj}).
\ep

\bL{L:Fclosed, F'not}  Let $\A$ be a simple, $\sigma$-unital, non-unital,   C*-algebra,  with metrizable $\TA$, and with projection-surjectivity and injectivity. Assume  that $F$ is a closed face, $\Sc(\tau)=\infty$ for all $\tau\in F$,  and the complementary face $F'$ is not closed. Then $\Icon \ne \Ib$. 
\eL

\bp
Let $0< \gamma < \min \Sc$. Then by Corollary  \ref {C:complementation} (ii), the function
$\frac{\gamma}{2}\mid_F\overset{\cdot}{+}\gamma\mid_{F'}$ belongs to $\LAff_{++}$ and is complemented under $\Sc$. Therefore there is a projection $P\not \in \A$ such that $\hat P(\tau)=\begin{cases}\frac{\gamma}{2}&\tau\in F\\\gamma&\tau\in F'\end{cases}$. Since $\hat P(\tau) \le \gamma$ for all $\tau$, $P\in \Ib$. Notice that $\hat P(\tau) < \gamma$ for every $\tau\not\in F'$. Since $F'$ is not closed, $\hat P$ is not continuous and hence $P\not \in \Icon$ by (\ref {e:cont proj}).

\ep
If the scale of $\A$ is not quasicontinuous and $\A$ is projection-surjective and injective, then that at least one of the inclusions $\Icon\subset \Ib\subset \Ifin$ must be proper. 
\bP{P: QCS} Let $\A$ be a simple, $\sigma$-unital, non-unital,  C*-algebra,  with metrizable tracial simplex, and projection-surjectivity and injectivity. Then
\item[(i)] If $F_\infty$ is finite and $F_\infty'$ is not closed then  $\Icon \ne \Ib$.
\item[(ii)] If $F_\infty$ is finite, $F_\infty'$ is closed, and $\Sc\mid _{F_\infty'}$ is not continuous, then $\Ifin\ne \Icon$. 
\item[(iii)] If $F_\infty$ is infinite and countable, then $I_b\ne \Ifin$.
\item[(iv)] If $F_\infty$ is uncountable, then $\Icon \ne \Ib$.\\
Thus if  $\Icon=\Ifin$, then $\A$ has quasicontinuous scale.

\eP
\bp 
\item [(i)] If $F_\infty$ is finite, then $F= \co(F_\infty)$ is closed and the conclusion is given by Lemma \ref {L:Fclosed, F'not}.
\item [(ii)] By Corollary \ref {C:complementation} the function $1\mid_{F_\infty}\overset{\cdot}{+} \frac{\Sc}{2}\in \LAff_{++}$ is complemented under $\Sc$ and therefore there is a projection $P$ such that $\hat P= 1\mid_{F_\infty}\overset{\cdot}{+} \frac{\Sc}{2}$. As $$\hat P(\tau)=\begin {cases} 1& \tau \in F_\infty\\ \frac{\Sc(\tau)}{2}< \infty & \tau\in F_\infty' \cap \Ext, \end{cases} $$
we see that $P\in \Ifin$. However $\hat P= \frac{\Sc}{2}$ on $F_\infty' $ is not continuous, and hence $P\not \in \Icon$.
\item [(iii)] Let $F_\infty= \{\tau_n\}$ and apply Lemma \ref {L: chain below} to the function $h=\Sc$ and the sequence $ \{\tau_n\}= F_\infty$. Then $\Sc=G+F$ where $G$ and $F$ are in $\LAff_{++}$. $G$ being complemented under $\Sc$, there is a projection $P\not \in \A$ such that $\hat P= G$. Then $P\in \Ifin$ because $G(\tau_n)< \infty$ for all $n$, but $P\not \in \Ib$ because $\hat P(\tau_n)$ is unbounded. 
\item [(iv)] By the assumption that $\TA$ is metrizable, we can find an element $x\in F_\infty$ that belongs to the closure of $ F_\infty\setminus \{x\}$.  Then $F:=\{x\}$ is closed, but  $F'\supset (F_\infty\setminus \{x\})$ is not closed, hence the conclusion follows again from Lemma \ref {L:Fclosed, F'not}. 
 \ep
 
 Notice that the proof of  (i) and (ii)  did not require metrizability.

We can sharpen the result of Proposition \ref {P: QCS} in the case when $\A$ is stable and hence $F_\infty = \Ext$. Then $\A$ has quasicontinuous scale if and only if $\Ext$ is finite.

\bP{P: stable} Let $\A$ be a simple, $\sigma$-unital,  C*-algebra,  with metrizable tracial simplex, and projection-surjectivity and injectivity and assume that $F_\infty = \Ext$. 
\item [(i)] $\Ext$ is finite if and only if $\Icon = \Ib$.
\item [(ii)] If furthermore $\TA$ is a Bauer simplex, then $\Ext$ is finite if and only if $\Ib = \Ifin$
\eP
\bp 
The necessity in both cases is given by Corollary \ref {C: fin bdry}.
\item [(i)] 
For the sufficiency, by  Proposition \ref {P: QCS} (iv), it is enough to prove that if $\Ext$ is infinite and countable then $\Icon \ne  \Ib$. To obtain that it is sufficient (and by Lemma \ref {L: Ifin not Ib} also necessary) to find a projection $P\in \Ib\setminus \Icon$.  By the surjectivity of $\A$ and the fact that every function in $ \LAff_{++}$ is complemented under  $\Sc$ because $\Sc(\tau)=\infty$ for all $\tau$, by  Corollary \ref {C: within Ib} and (\ref {e:cont proj})  we just need to  construct a bounded function $g\in \LAff_{++}\setminus \Aff$. 

Let $\{\tau_j\}_1^\infty$ be an enumeration of $\Ext$. Let $X_n:=\co\{\tau_1, \cdots, \tau_n\}$, then $X_n$ is closed and hence it is a split face. Define
$$f_n:= 1\mid_{X_n}\+2\mid_{X_n'}$$
By Lemma \ref {L:direct sum}, $f_n\in\LAff_{++}$, $f_n \le 2$ and clearly, $f_n$ is monotone noincreasing. If for some $n$, the function $f_n$ is not continuous, then we are done. Assume therefore that all the functions $f_n$ are continuous and let $f:=\lim_nf_n$. We claim that $f$ is not continuous. Indeed by the compactness of $\TA$, there is a subsequence $\tau_{j_k} $ that converges to some $\mu\in \TA$. Then for every $n$, 
$f_n(\tau_{j_k})\to f_n(\mu)$. Since $f_n(\tau_{j_k})=2$ for $j_k\ge n$ , we thus have $f_n(\mu)=2$ and hence $f(\mu)=2$. On the other hand, for every $k$, 
$$f(\tau_{j_k})= \lim _nf_n(\tau_{j_k})=1.$$
As a consequence the function $g:= 3-f$ is bounded but also not continuous. Since $g= \lim_n(3-f_n)$ is an increasing limit of functions in $\Aff$, it follows that $g \in \LAff_{++}$ which concludes the proof.
\item [(ii)]  Reasoning as in part (i), it is enough to assume that $\Ext$ is infinite (and uncountable) and then construct a function $f\in \LAff_{++}$ which is finite but unbounded on $\Ext$. We start by choosing a strictly positive lower semicontinuous function  $\tilde f : \Ext \mapsto (0, \infty)$  on the compact set $\Ext$ which is finite and unbounded.  
For instance,  let $d$ be the metric of $\TA$ restricted to $\Ext$,  $\tau_o \in \Ext$ be an accumulation point of $\Ext$ and set
$$
\tilde f(\tau) := 
\begin{cases}
1/d(\tau,\tau_o) & \tau \neq \tau_o \\
1 & \tau = \tau_o.
\end{cases}
$$
It is easily seen that $\tilde f$ satisfies the required conditions. Decompose 
$\tilde f = \sum_{n=1}^{\infty} \tilde f_n$ as a pointwise converging sum of functions  $\tilde f_n \in \mathrm{Aff}(\Ext)_{++} $ (Proposition \ref {P: sup aff}) .
By \cite[Corollary 11.15] {Goodearl}, for each $n$, there is a $f_n \in \Aff$
such that $f_n \big|_{\Ext} = \tilde f_n$,  and it is easy to see that $f_n$ must be strictly positive.
Then  $f : =\sum_{n=1}^{\infty} f_n \in \LAff_{++}$ and  $f\big |_{\Ext} = \tilde f$. Thus $f$ satisfies the required conditions.

\ep

If $\Icon\ne \Ifin$ we can draw several conclusions about $\Ma$.
\bL{L:above discont} Let $\A$ be a simple, $\sigma$-unital, non-unital,   C*-algebra,  with metrizable $\TA$, and with projection-surjectivity and injectivity. 
Assume $P$ is a projection in $\Ib\setminus \Icon$.  Then there is a projection $Q\in \Ib\setminus  \Icon$ such that
\item [(i)] $\hat P(\tau)< \hat Q(\tau)$ for every $\tau\in \TA.$
\item [(ii)] $I(P)= I(Q).$
\item [(iii)] $P\not\preceq Q.$\\
In particular, strict comparison of projections by traces does not hold on $\Ma$. 
\eL 

\bp
Let $P\in \Ib\setminus \Icon$. By (\ref {e:cont proj}) and Lemma  \ref {L: char Ifin, Ib}, $\hat P$ is bounded but not continuous. By invoking Proposition \ref {P:dividing} and recalling that if we divide $P$ into $n$ equivalent projection summands, each summand generates the same ideal as $P$, we can assume without loss of generality that $\sup \hat P < \min \Sc$.  Let $\sup \hat P < c < \min \Sc$. Then $\frac{1}{2}(\hat P+c)\in \LAff_{++}\setminus \Aff$.  Moreover, both $\Sc-c \in \LAff_{++}$ and $\Sc-\hat P= \widehat {1_{\Ma}-P}\in \LAff_{+}$. Since
$$\Sc= \frac{1}{2}(\hat P+c) +  \frac{1}{2}\big(\Sc-\hat P + \Sc- c\big),$$
$\frac{1}{2}(\hat P+c)$ is complemented under $\Sc$ and hence
 by the 1-projection-surjectivity of $\A$ there is a projection $Q$ such that $\hat Q= \frac{1}{2}(\hat P+c)$.  Again by (\ref {e:cont proj}) and Lemma  \ref {L: char Ifin, Ib}, $Q\in \Ib \setminus  \Icon$.  Condition (i) holds as  $\hat Q(\tau)-\hat P(\tau) =  \frac{1}{2}(c-\hat P(\tau))>0$ for all $\tau$. As $\hat P +c= 2\hat Q$,  it follows by Proposition \ref {P: in ideal} that  $P\in I(Q)$. Furthermore,   let $m\in \mathbb N$ be such that $(2m-1) \inf \hat P> c$. Then $g:= \frac{1}{2}\big ((2m-1)\hat P-c\big)\in \LAff_{++}$ and  $\hat Q+g = m\hat P$. Thus $Q\in I(P)$ by Proposition \ref {P: in ideal} and hence $I(P)=I(Q)$, which establishes condition (ii).

To prove (iii) assume by contradiction that $P\preceq Q$. Then there is a function  $f\in \LAff_{++}$ such that $\hat P+f= \hat Q$. But then $f= \frac{1}{2}(c-\hat P)$ by the boundedness of $\hat P$,  whence $f$ is also upper semicontinuous and hence it is continuous. This implies that $\hat P$ is continuous, a contradiction.
\ep

\bL{L: chains}  Let $\A$ be a simple, $\sigma$-unital, non-unital,  C*-algebra,  with metrizable $\TA$, and with projection-surjectivity and injectivity. Assume that there is a projection $P\in \Ib\setminus \Icon$ (resp.,  $P\in \Ifin\setminus \Ib$). Then there is a projection $P_1\in \Ib\setminus \Icon$ (resp., $P_1\in\Ifin\setminus \Ib$) such that $I(P_1)\subsetneq I(P)$. Therefore $I_b$ (resp. $\Ifin$), contains an infinite decreasing chain of principal ideals. \eL
\bp
Assume first that $P\in \Ib\setminus \Icon$.
By(\ref {e:cont proj}) and Lemma  \ref {L: char Ifin, Ib}, $\hat P$ is a bounded function in $ \LAff_{++}$ and it has at least one point of discontinuity $\mu\in \TA$.  Then by Lemma \ref {L: under discont}, 
 $\hat P = G+F$ where $G,F\in \LAff_{++}$ are both discontinuous at $\mu$ but for which there is a sequence $\tau_n\to \mu$ such that $G(\tau_n)\to G(\mu)$, and  $\hat P(\tau_n)\not \to \hat P(\mu)$. By Lemma  \ref {L:subequiv},  there is a projection $P_1$ such that $\widehat {P_1}= G$ and  $P_1\preceq P$. Then $P_1\in\Ib\setminus \Icon$ and   $I(P_1)\subset I(P)$. If 
$I(P_1)= I(P)$, we would have $P\in I(P_1)$ and hence by Proposition \ref {P: in ideal}  there would be an $m\in \mathbb N$ and a function $f\in \LAff_{++}$ such that  $\hat P + f= m\widehat {P_1}= m G$. However, since $mG(\tau_n)\to mG(\mu)$ and both  $\hat P$ and $  f$ are lower semicontinuous, we would conclude by Lemma  \ref {L:facts} that $\hat P(\tau_n)\to \hat P(\mu)$, a contradiction. 

Assume now that $P\in \Ifin\setminus \Ib$.
By Lemma \ref {L: chain below}  there is sequence $\tau_n \in \Ext$ such that $\hat P(\tau_n)$ is finite for every $n$ but the sequence is unbounded. Apply Lemma \ref {L: chain below} to the function $h:=\hat P$ and the sequence $\{\tau_n\}$ to decompose $h= G+F$ into the sum of $G,F\in \LAff_{++}$, with $G$ unbounded but $\sup_n\frac{\hat P(\tau_n)}{G(\tau_n)}=\infty$.  Then there is a  projection $P_1 \preceq P$ with $\widehat {P_1}= G$ and hence 
$P_1\in \Ifin\setminus \Ib$. Furthermore,  $P\not\in I( P_1)$. Indeed, otherwise there would be an $g\in \LAff_{++}$ and $m\in \mathbb N$ such that  $\hat P+f=mG$. But then $\frac{\hat P(\tau_n)}{G(\tau_n)}\le m $ for every $n$, a contradiction.
\ep

\bC{C:fin many  ideals}
Let $\A$ be a simple, $\sigma$-unital, non-unital, C*-algebra,  with metrizable $\TA$, and with projection-surjectivity and injectivity. If $\Icon\ne\Ifin$ then $\Ma$ has infinitely many (principal) ideals and therefore $V(\Ma)$ contains infinitely many (principal) order ideals.\eC
\bp
If $\Icon\ne\Ifin$, then at least one of the inclusions $\Icon\subset \Ib\subset \Ifin$ must be proper. By Lemma \ref {L: Ifin not Ib}, there must be a projection in $\Ifin\setminus \Ib$ or in $\Ib\setminus \Icon$. In either case the conclusion follows from Lemma \ref {L: chains}. By Lemma \ref {L: id vs oid} we see that $V(\Ma)$ contains infinitely many (principal) order ideals.
\ep

Notice that the chains of principal ideal constructed in Lemma \ref {L: chains}  are decreasing. If $\A$ is stable and has countably infinite extremal boundary, we can also construct increasing chains.

\bP{P:ext count} Let $\A$ be  a simple, stable, $\sigma$-unital, non-unital, C*-algebra,  with metrizable  tracial simplex $\TA$, countably infinite extremal boundary $\Ext$,  and projection-surjectivity and injectivity.
For every projection $P\in \Ifin$ there is a continuous chain of projections $P_t\in  \Ifin\setminus \Ib$ for $t\ge1$ such that
$$I(P)\subsetneq I(P_1) \subsetneq I(P_s)\subsetneq I(P_t)\quad\forall\,~ 1< s< t.$$ 
\eP
\bp
 Let $\{\tau_j\}_1^\infty$ be an enumeration of $\Ext$. Since $0< \hat P(\tau_n)< \infty$, we can find a sequence $\beta_n$ such that
 \begin{itemize} \item [(i)] $1< (\beta_n-1)\hat P(\tau_n)$ is monotone nondecreasing
  \item [(ii)] $ \beta_n\hat P(\tau_n)$ is monotone nondecreasing
  \item  [(iii)] $\beta_n\to \infty$
\end{itemize} 
By Corollary \ref {C:increasing sequences} there exist a projections $P_1$ and a function $g\in \LAff_{++}$  such that  for every $n$ \ba \hat P_1(\tau_n)&= \beta_n\hat P(\tau_n)\\  
g(\tau_n)&= (\beta_n-1)\hat P(\tau_n).\end{align*} Since $\hat P(\tau_n)+ g (\tau_n)= \hat P_1(\tau_n)$, it follows by  Proposition \ref{P: sup aff} (iii) that 
$\hat P+ g= \hat P_1$. Then $I(P)\subset I(P_1) $ by Lemma  \ref {L:subequiv} and Proposition \ref {P: in ideal}. On the other hand $\sup_n\frac{\hat P_1(\tau_n)}{\hat P(\tau_n)}=\infty$ and hence by Proposition \ref {P: in ideal} it follows that  $P_1\not\in I(P)$.

Next, for every $t>1$, let $P_t$ be the projection for which $\hat P_t(\tau_n):= \big ( \beta_n\hat P(\tau_n)\big )^t$. 
Since  for $1\le s< t< \infty$ the sequence 
$$\big ( \beta_n\hat P(\tau_n)\big )^t-\big ( \beta_n\hat P(\tau_n)\big )^s= \big( \beta_n\hat P(\tau_n)\big )^s\Big (\big( \beta_n\hat P(\tau_n)\big )^{t-s}- 1\Big)$$
is monotone nondecreasing, again by Corollary \ref {C:increasing sequences},   there exists a function $g\in \LAff_{++}$ that achieves the values of that sequence at $\tau_n$, that is
$$
\hat P_s(\tau_n)+g(\tau_n)= \hat P_t(\tau_n).
$$
But then, again by Proposition \ref{P: sup aff} (iii) and Proposition \ref {P: in ideal} it follows that $\hat P_s+g=\hat P_t$, hence $P_s\preceq P_t$ and thus $I(P_s)\subset I(P_t)$. Since $\sup_n\frac{\hat P_t(\tau_n)}{\hat P_s(\tau_n)}=\infty$, again  by Proposition \ref {P: in ideal} it follows that $P_t\not\in I(P_s)$.
\ep

We collect now the results obtained in this section in  our main theorem.  
\bT{T:main}
Let $\A$ be a simple, $\sigma-$unital, non-unital  C*-algebra, with metrizable tracial simplex $\TA$, projection-surjectivity and injectivity, and strict comparison of positive elements by traces. Then the following are equivalent 
\item [(i)] $\A$ has quasicontinuous scale;
\item [(ii)] $\Ma$ has strict comparison of positive elements by traces;
\item [(iii)] $\Ma/\A$ is purely infinite;
\item [(iii$'$)] $\Ma/\Imin$ is purely infinite;
\item [(iv)] $\Ma$ has finitely many ideals;
\item [(v)] $\Imin=\Ifin$.\\
Consider in addition 
\item [(vi)] $V(\Ma)$ has finitely many order ideals.
Then (vi) implies (i)-(v). If $\A$ is n-projection-surjective and n-projection-injective for every $n$, then (vi) is equivalent to  (i)-(v).
\eT
We will always assume that $\A$ is simple, $\sigma$-unital but not unital, non-elementary, and with non-empty tracial simplex (and hence stably finite), but not all of the other three  hypotheses (metrizability of $\TA$, projection-surjectivity and injectivity of $\A$, and strict comparison of positive elements of $\A$), will be necessary for all the implications.  In the proofs of the various implications, we will list which of these other hypotheses are used and/or which ones can be weakened.

\bp 

 (i) $\Rightarrow$ (ii). By Theorem \ref {T:strict compar} (\cite[Theorem 6.6]{KNZCompPos}). For this implication we need only strict comparison of positive elements for $\A$.\\
 (ii) $\Rightarrow$ (iii). By Proposition \ref {P:strict comp implies prop inf}. For this implication we do not require any of the other three hypotheses.\\
 (iii) $\Leftrightarrow$ (iii$'$). 
In view of the exact sequence
$$ 0\to \Imin/\A\to\Ma/\A\to \Ma/\Imin\to 0$$ the conclusion follows  from the ``two out of three" property  (\cite [Theorem 4.19]{KirchRor}) provided that   $\Imin/\A$ is purely infinite. By \cite [Theorem 4.8]{KNZMin}), a sufficient condition for  $\Imin/\A$ to be purely infinite is that $\A$ is non-elementary and that $\Imin\ne \A$, which follows from the strict comparison of positive elements in $\A$ (\cite [Corollary 3.15, Proposition 5.4, Theorem 5.6 ]{KNZMin}). If $\A$ is separable then also $\Imin\ne \A$ (\cite [Corollary 3.15]{KNZMin}) so we can replace the condition of strict comparison of positive elements in $\A$ with the separability of $\A$.\\
(iii) $\Rightarrow$ (v) 
Proposition \ref {P: corona pi}, which does not require any of three additional hypotheses, guarantees that if $\Ma/\A$ is purely infinite, then all the projections of  $\Ifin$ belong to $\Icon$. This in turn implies that $\Icon=\Ifin$ by Lemma \ref {L: Ifin not Ib},  which makes use only 1-projection-surjectivity for $\A$.  Finally, $\Imin=\Icon$ by \cite [Theorem 5.6] {KNZMin}) which depends only on strict comparison of positive elements.\\
(v)$\Rightarrow$ (i) By Proposition \ref {P: QCS}. Projection-surjectivity and injectivity for $\A$ and metrizability of $\TA$ are used for obtaining that $\Icon=\Ifin$ implies quasicontinuity of the scale. As above, strict comparison of positive elements is used for obtaining that $\Imin=\Icon$ (\cite [Theorem 5.6] {KNZMin}).\\
(i)$\Rightarrow$ (iv) By Corollary \ref{C: structure ideals}, which makes use only of strict comparison of positive elements for $\A$.\\
(iv)$\Rightarrow$ (v). By Corollary \ref {C:fin many  ideals}. For this implication we 
use projection-surjectivity and injectivity and the metrizability of $\TA$ to obtain that $\Icon =\Ifin$, and again strict comparison of positive elements to obtain that $\Imin=\Icon$.\\ 
(vi)$\Rightarrow$ (v). Strict comparison on $\A$ guarantee that $\Imin=\Icon$ and  metrizability of $\TA$ and projection-surjectivity and injectivity permit to apply Lemma \ref {L: Ifin not Ib} and Lemma \ref {L: chains}. Thus if $\Imin\ne \Ifin$ then  $\Ifin$ contains an infinite chain of {\it principal} ideals and hence $V(\Ma)$ has an infinite chain of (principal) order ideals. \\
If $\A$ is n-projection-surjective and n-projection-injective for every $n$, then (i) $\Rightarrow$ (vi) by Corollary \ref {C:order ideals} which requires metrizability of $\TA$. 

\ep 

When the algebra $\A$ is separable and stable, asking for strict comparison is redundant (Theorem \ref {T: proj surj inj implies strict comp}) and we see that $\Imin=\Ifin$ if and only if $\Imin= \Ib$.
\bC{C:main}
Let $\A$ be a simple, separable, stable,  C*-algebra, with projection-surjectivity and injectivity. Then the following are equivalent 
\item [(i)] The extremal boundary $\Ext$ is finite;
\item [(ii)] $\Ma$ has strict comparison of positive elements by traces;
\item [(iii)] $\Ma/\A$ is purely infinite;
\item [(iii$'$)] $\Ma/\Imin$ is purely infinite;
\item [(iv)] $\Ma$ has finitely many ideals;
\item [(v)] $\Imin=\Ifin$;
\item [(v$'$)] $\Imin= \Ib$.
\item [(vi)] $V(\Ma)$ has finitely many order ideals;
\eC
\bp 
All the hypotheses of Theorem \ref {T:main} are satisfied: metrizability is implied by the separability of $\A$,  strict comparison of positive elements for $\A$ is implied 
projection-surjectivity and injectivity (Theorem \ref {T: proj surj inj implies strict comp}).  Thus conditions (i), (ii),(iii), (iii$'$), (iv), (v), and (vi) are equivalent, where for (i) we notice that for stable C*-algebras quasicontinuity of the scale is equivalent to finiteness of the extremal boundary. \\
 (i) $\Leftrightarrow$ (v$'$) By Proposition \ref {P: stable}.
 
\ep

\bR{R:KP-KNP} For the class of algebras $\A$ that are separable, non-unital, simple, have real rank zero, stable rank one, and have weakly unperforated $K_0(\A)$, 
the equivalence  of (i), (iii), (iii$'$), and (v) was established in \cite[Theorem 3.4]{KucPer} under the additional condition that $\A$ has finitely many infinite extremal quasitraces.
For the same class of algebras, the equivalence  of the above conditions  with (iv) was established in \cite[Theorem 3.6]{KucPer} under the additional condition that $\A$ is exact, is the stabilization of a unital algebra, and $\TA$ is a Bauer simplex.
In \cite{KucNgPer} The equivalence  of (i), (iii), (iv) was established under the condition that $\A$ is the stabilization of a unital algebra, is separable, simple, and is either exact and $\Z$-stable or an AH-algebra with slow dimension growth. These results, in turn, are generalizations of earlier work in \cite {ElliottMatroidII}, \cite {LinMultIdeals}, \cite {LinSimpl}, \cite {LinContScale}, \cite {LinSimple}, \cite {RordamIdeals}, \cite {SZ90}, \cite {ZhangWeylVonNeumann}
\eR

\end{document}